\catcode`\@=11

\font\tenmsa=msam10
\font\sevenmsa=msam7
\font\fivemsa=msam5
\font\tenmsb=msbm10
\font\sevenmsb=msbm7
\font\fivemsb=msbm5
\newfam\msafam
\newfam\msbfam
\textfont\msafam=\tenmsa  \scriptfont\msafam=\sevenmsa
  \scriptscriptfont\msafam=\fivemsa
\textfont\msbfam=\tenmsb  \scriptfont\msbfam=\sevenmsb
  \scriptscriptfont\msbfam=\fivemsb

\def\hexnumber@#1{\ifnum#1<10 \number#1\else
 \ifnum#1=10 A\else\ifnum#1=11 B\else\ifnum#1=12 C\else
 \ifnum#1=13 D\else\ifnum#1=14 E\else\ifnum#1=15 F\fi\fi\fi\fi\fi\fi\fi}

\def\msa@{\hexnumber@\msafam}
\def\msb@{\hexnumber@\msbfam}
\mathchardef\boxdot="2\msa@00
\mathchardef\boxplus="2\msa@01
\mathchardef\boxtimes="2\msa@02
\mathchardef\square="0\msa@03
\mathchardef\blacksquare="0\msa@04
\mathchardef\centerdot="2\msa@05
\mathchardef\lozenge="0\msa@06
\mathchardef\blacklozenge="0\msa@07
\mathchardef\circlearrowright="3\msa@08
\mathchardef\circlearrowleft="3\msa@09
\mathchardef\rightleftharpoons="3\msa@0A
\mathchardef\leftrightharpoons="3\msa@0B
\mathchardef\boxminus="2\msa@0C
\mathchardef\Vdash="3\msa@0D
\mathchardef\Vvdash="3\msa@0E
\mathchardef\vDash="3\msa@0F
\mathchardef\twoheadrightarrow="3\msa@10
\mathchardef\twoheadleftarrow="3\msa@11
\mathchardef\leftleftarrows="3\msa@12
\mathchardef\rightrightarrows="3\msa@13
\mathchardef\upuparrows="3\msa@14
\mathchardef\downdownarrows="3\msa@15
\mathchardef\upharpoonright="3\msa@16

\mathchardef\downharpoonright="3\msa@17
\mathchardef\upharpoonleft="3\msa@18
\mathchardef\downharpoonleft="3\msa@19
\mathchardef\rightarrowtail="3\msa@1A
\mathchardef\leftarrowtail="3\msa@1B
\mathchardef\leftrightarrows="3\msa@1C
\mathchardef\rightleftarrows="3\msa@1D
\mathchardef\Lsh="3\msa@1E
\mathchardef\Rsh="3\msa@1F
\mathchardef\rightsquigarrow="3\msa@20
\mathchardef\leftrightsquigarrow="3\msa@21
\mathchardef\looparrowleft="3\msa@22
\mathchardef\looparrowright="3\msa@23
\mathchardef\circeq="3\msa@24
\mathchardef\succsim="3\msa@25
\mathchardef\gtrsim="3\msa@26
\mathchardef\gtrapprox="3\msa@27
\mathchardef\multimap="3\msa@28
\mathchardef\therefore="3\msa@29
\mathchardef\because="3\msa@2A
\mathchardef\doteqdot="3\msa@2B

\mathchardef\triangleq="3\msa@2C
\mathchardef\precsim="3\msa@2D
\mathchardef\lesssim="3\msa@2E
\mathchardef\lessapprox="3\msa@2F
\mathchardef\eqslantless="3\msa@30
\mathchardef\eqslantgtr="3\msa@31
\mathchardef\curlyeqprec="3\msa@32
\mathchardef\curlyeqsucc="3\msa@33
\mathchardef\preccurlyeq="3\msa@34
\mathchardef\leqq="3\msa@35
\mathchardef\leqslant="3\msa@36
\mathchardef\lessgtr="3\msa@37
\mathchardef\backprime="0\msa@38
\mathchardef\risingdotseq="3\msa@3A
\mathchardef\fallingdotseq="3\msa@3B
\mathchardef\succcurlyeq="3\msa@3C
\mathchardef\geqq="3\msa@3D
\mathchardef\geqslant="3\msa@3E
\mathchardef\gtrless="3\msa@3F
\mathchardef\sqsubset="3\msa@40
\mathchardef\sqsupset="3\msa@41
\mathchardef\trianglerighteq="3\msa@44
\mathchardef\trianglelefteq="3\msa@45
\mathchardef\bigstar="0\msa@46
\mathchardef\between="3\msa@47
\mathchardef\blacktriangledown="0\msa@48
\mathchardef\blacktriangleright="3\msa@49
\mathchardef\blacktriangleleft="3\msa@4A
\mathchardef\blacktriangle="0\msa@4E
\mathchardef\triangledown="0\msa@4F
\mathchardef\eqcirc="3\msa@50
\mathchardef\lesseqgtr="3\msa@51
\mathchardef\gtreqless="3\msa@52
\mathchardef\lesseqqgtr="3\msa@53
\mathchardef\gtreqqless="3\msa@54
\mathchardef\Rrightarrow="3\msa@56
\mathchardef\Lleftarrow="3\msa@57
\mathchardef\veebar="2\msa@59
\mathchardef\barwedge="2\msa@5A
\mathchardef\doublebarwedge="2\msa@5B
\mathchardef\angle="0\msa@5C
\mathchardef\measuredangle="0\msa@5D
\mathchardef\sphericalangle="0\msa@5E
\mathchardef\varpropto="3\msa@5F
\mathchardef\smallsmile="3\msa@60
\mathchardef\smallfrown="3\msa@61
\mathchardef\Subset="3\msa@62
\mathchardef\Supset="3\msa@63
\mathchardef\Cup="2\msa@64

\mathchardef\Cap="2\msa@65

\mathchardef\curlywedge="2\msa@66
\mathchardef\curlyvee="2\msa@67
\mathchardef\leftthreetimes="2\msa@68
\mathchardef\rightthreetimes="2\msa@69
\mathchardef\subseteqq="3\msa@6A
\mathchardef\supseteqq="3\msa@6B
\mathchardef\bumpeq="3\msa@6C
\mathchardef\Bumpeq="3\msa@6D
\mathchardef\lll="3\msa@6E

\mathchardef\ggg="3\msa@6F

\mathchardef\circledS="0\msa@73
\mathchardef\pitchfork="3\msa@74
\mathchardef\dotplus="2\msa@75
\mathchardef\backsim="3\msa@76
\mathchardef\backsimeq="3\msa@77
\mathchardef\complement="0\msa@7B
\mathchardef\intercal="2\msa@7C
\mathchardef\circledcirc="2\msa@7D
\mathchardef\circledast="2\msa@7E
\mathchardef\circleddash="2\msa@7F
\def\ulcorner{\delimiter"4\msa@70\msa@70 }
\def\urcorner{\delimiter"5\msa@71\msa@71 }
\def\llcorner{\delimiter"4\msa@78\msa@78 }
\def\lrcorner{\delimiter"5\msa@79\msa@79 }
\def\yen{\mathhexbox\msa@55 }
\def\checkmark{\mathhexbox\msa@58 }
\def\circledR{\mathhexbox\msa@72 }
\def\maltese{\mathhexbox\msa@7A }
\mathchardef\lvertneqq="3\msb@00
\mathchardef\gvertneqq="3\msb@01
\mathchardef\nleq="3\msb@02
\mathchardef\ngeq="3\msb@03
\mathchardef\nless="3\msb@04
\mathchardef\ngtr="3\msb@05
\mathchardef\nprec="3\msb@06
\mathchardef\nsucc="3\msb@07
\mathchardef\lneqq="3\msb@08
\mathchardef\gneqq="3\msb@09
\mathchardef\nleqslant="3\msb@0A
\mathchardef\ngeqslant="3\msb@0B
\mathchardef\lneq="3\msb@0C
\mathchardef\gneq="3\msb@0D
\mathchardef\npreceq="3\msb@0E
\mathchardef\nsucceq="3\msb@0F
\mathchardef\precnsim="3\msb@10
\mathchardef\succnsim="3\msb@11
\mathchardef\lnsim="3\msb@12
\mathchardef\gnsim="3\msb@13
\mathchardef\nleqq="3\msb@14
\mathchardef\ngeqq="3\msb@15
\mathchardef\precneqq="3\msb@16
\mathchardef\succneqq="3\msb@17
\mathchardef\precnapprox="3\msb@18
\mathchardef\succnapprox="3\msb@19
\mathchardef\lnapprox="3\msb@1A
\mathchardef\gnapprox="3\msb@1B
\mathchardef\nsim="3\msb@1C
\mathchardef\napprox="3\msb@1D
\mathchardef\nsubseteqq="3\msb@22
\mathchardef\nsupseteqq="3\msb@23
\mathchardef\subsetneqq="3\msb@24
\mathchardef\supsetneqq="3\msb@25
\mathchardef\subsetneq="3\msb@28
\mathchardef\supsetneq="3\msb@29
\mathchardef\nsubseteq="3\msb@2A
\mathchardef\nsupseteq="3\msb@2B
\mathchardef\nparallel="3\msb@2C
\mathchardef\nmid="3\msb@2D
\mathchardef\nshortmid="3\msb@2E
\mathchardef\nshortparallel="3\msb@2F
\mathchardef\nvdash="3\msb@30
\mathchardef\nVdash="3\msb@31
\mathchardef\nvDash="3\msb@32
\mathchardef\nVDash="3\msb@33
\mathchardef\ntrianglerighteq="3\msb@34
\mathchardef\ntrianglelefteq="3\msb@35
\mathchardef\ntriangleleft="3\msb@36
\mathchardef\ntriangleright="3\msb@37
\mathchardef\nleftarrow="3\msb@38
\mathchardef\nrightarrow="3\msb@39
\mathchardef\nLeftarrow="3\msb@3A
\mathchardef\nRightarrow="3\msb@3B
\mathchardef\nLeftrightarrow="3\msb@3C
\mathchardef\nleftrightarrow="3\msb@3D
\mathchardef\divideontimes="2\msb@3E
\mathchardef\varnothing="0\msb@3F
\mathchardef\nexists="0\msb@40
\mathchardef\mho="0\msb@66
\mathchardef\thorn="0\msb@67
\mathchardef\beth="0\msb@69
\mathchardef\gimel="0\msb@6A
\mathchardef\daleth="0\msb@6B
\mathchardef\lessdot="3\msb@6C
\mathchardef\gtrdot="3\msb@6D
\mathchardef\ltimes="2\msb@6E
\mathchardef\rtimes="2\msb@6F
\mathchardef\shortmid="3\msb@70
\mathchardef\shortparallel="3\msb@71
\mathchardef\smallsetminus="2\msb@72
\mathchardef\thicksim="3\msb@73
\mathchardef\thickapprox="3\msb@74
\mathchardef\approxeq="3\msb@75
\mathchardef\succapprox="3\msb@76
\mathchardef\precapprox="3\msb@77
\mathchardef\curvearrowleft="3\msb@78
\mathchardef\curvearrowright="3\msb@79
\mathchardef\digamma="0\msb@7A
\mathchardef\varkappa="0\msb@7B
\mathchardef\hslash="0\msb@7D
\mathchardef\hbar="0\msb@7E
\mathchardef\backepsilon="3\msb@7F
\def\Bbb{\ifmmode\let\next\Bbb@\else
 \def\next{\errmessage{Use \string\Bbb\space only in math mode}}\fi\next}
\def\Bbb@#1{{\Bbb@@{#1}}}
\def\Bbb@@#1{\fam\msbfam#1}

\catcode`\@=\active

\def\del{\partial}

\def\CR{\hbox{{$\cal R$}}}

\def\cu{{\upsilon}} 
\def\cv{{\vartheta}} 

\def\lform{\hbox{$\sqcup$}\llap{\hbox{$\sqcap$}}}

\def\h{{{1\over2}}}

\def\R{{\Bbb R}}
\def\C{{\Bbb C}}
\def\Z{{\Bbb Z}}

\def\eps{{\epsilon}}

\def\dcross{{\bowtie}}
\def\codcross{{\blacktriangleright\!\!\blacktriangleleft}}
\def\rbiprod{{\cdot\kern-.33em\triangleright\!\!\!<}}
\def\lbiprod{{>\!\!\!\triangleleft\kern-.33em\cdot\, }}

\def\tens{\mathop{\otimes}}

\def\la{{\triangleright}}

\def\isom{{\cong}}

\def\ev{{\rm ev}}

\def\id{{\rm id}}

\def\<{\langle}
\def\>{\rangle}
\def\dila{{\varsigma}}

\def\equad{\kern -1.7em}

\def\eqn#1#2{\begin{equation}#2\label{#1}\end{equation}}

\def\haj#1{{\mathaccent20 {#1}}}
\def\Vhaj{{V\haj{\ }}}

\def\o{{}_{\scriptscriptstyle(1)}}
\def\t{{}_{\scriptscriptstyle(2)}}

\def\bo{{}^{\bar{\scriptscriptstyle(1)}}}
\def\bt{{}^{\bar{\scriptscriptstyle(2)}}}

\def\und#1{{\underline {#1}}}

\def\uo{{{}^{\scriptscriptstyle(1)}}}
\def\ut{{{}^{\scriptscriptstyle(2)}}}

\def\umo{{{}^{\scriptscriptstyle-(1)}}}
\def\umt{{{}^{\scriptscriptstyle-(2)}}}

\def\Bo{{{}_{\und{\scriptscriptstyle(1)}}}}
\def\Bt{{{}_{\und{\scriptscriptstyle(2)}}}}

\def\text#1{\mbox{\rm #1}}
\def\note#1{}

\def\blacksquare{{\lform}}
\def\frac#1#2{{{#1\over#2}}}

\def\ceqn#1#2{\begin{equation}\label{#1}
\begin{array}{c}#2\end{array}\end{equation}}

\def\vect{{\bf t}}\def\vecs{{\bf s}}
\def\vecu{{\bf u}}\def\vecx{{\bf x}}\def\vecp{{\bf p}}
\def\vecl{{\bf l}}
\def\vecL{{\bf L}}
\def\vecm{{\bf m}}\def\vecP{{\bf P}}

\def\FRT{R\vect_1\vect_2=\vect_2\vect_1R}
\def\bramat{R_{21}\vecu_1 R\vecu_2=\vecu_2 R_{21}\vecu_1 R}
\def\eucmat{R_{21}\vecx_1 \vecx_2=\vecx_2\vecx_1 R}
\def\mulstat{R^{-1}\vecu'_1 R\vecu_2=\vecu_2 R^{-1}\vecu'_1 R}
\def\addstat{R^{-1}\vecu_1'R\vecu_2= \vecu_2 R_{21}\vecu_1'R}

\def\mink{{\R_q^{1,3}}}

\documentstyle[11pt]{article} \textheight 23.6cm \textwidth 16cm
\evensidemargin
0in \topskip 28pt

\newtheorem{lemma}{Lemma}[section] \newtheorem{propos}[lemma]{Proposition}
 \newtheorem{theorem}[lemma]{Theorem}
\newtheorem{cor}[lemma]{Corollary}

\begin{document}\baselineskip 20pt

{\ }\hskip 4.2in DAMTP/95-08 (revised)
\vspace{.2in}

\begin{center} {\LARGE SOME REMARKS ON THE q-POINCARE ALGEBRA}\\ {\LARGE IN
R-MATRIX FORM} \\ \baselineskip 13pt
{\ } {\ }\\ S.  Majid\footnote{Royal Society University Research Fellow and
Fellow of Pembroke College, Cambridge}\\{\ }\\ Department of Applied
Mathematics \& Theoretical Physics\\ University of Cambridge, Cambridge CB3
9EW\footnote{On leave calendar years 1995+1996 at  Department of Mathematics,
Harvard University}  \end{center}

\vspace{10pt} \begin{quote}\baselineskip 12pt
\noindent{\bf Abstract} The braided approach to q-deformation (due to the
author and collaborators) gives natural algebras  $\bramat$ and $\eucmat$ for
q-Minkowski and
q-Euclidean spaces respectively. These algebras are covariant under a
corresponding background `rotation' quantum group.  Semidirect
product by this according to the bosonisation procedure (also due to the
author) gives the corresponding Poincar\'e quantum groups.
We review the construction and collect the resulting R-matrix formulae for both
Euclidean and Minkowski cases in both enveloping algebra and function algebra
form, and the duality between them. Axioms for the Poincar\'e quantum group
$*$-structure and the dilaton problem are discussed.

\end{quote} \baselineskip 21.5pt

\section{Introduction}

The programme of q-deforming the basic geometrical notions of spacetime has
been extensively studied in
recent years and by several groups. Of the various approaches, two have been
pushed quite far. One, which
is the `minimal deformation', involves a non-commutative time co-ordinate but
the space co-ordinates remain
unchanged. The corresponding Poincar\'e quantum group is the so-called
`$\kappa$-Poincar\'e' studied by Lukierski {\em et al}\cite{LNRT:def} and
others,  obtained
by a contraction procedure. Its semidirect
product structure in terms of (usual) Lorentz rotations and a $\kappa$-deformed
4-momentum is due to the author and Ruegg\cite{MaRue:bic} and came later, along
with the identification of the correct (non-commutative) Minkowski space such
that the $\kappa$-Poincar\'e algebra acts covariantly on it. The semidirect
structure is an example of a {\em bicrossproduct Hopf algebra} as introduced
in\cite{Ma:phy}. We don't want to say too much about this here except that this
approach uses the standard theory of Hopf algebras (or quantum groups) and not
the more novel braided group theory.

The second deformation, which is the one that concerns us in this paper, is a
programme introduced first by Carow-Watamura {\em et al.}
\cite{CWSSW:lor}\cite{CWSSW:ten} from consideration of the tensor product of
two quantum planes (the 2-spinorial or twistor point of view). At about the
same time in 1990, the present author introduced independently a theory of
braided groups\cite{Ma:rec}\cite{Ma:bra}, including braided matrices
$B(R)$\cite{Ma:exa} with relations
\eqn{bramat}{\vecu=\{u^i{}_j\},\quad \bramat}
which for the $2\times 2$ case also provides a natural definition of
$q$-deformed Minkowski space. There is a hermitian $*$-structure whenever $R$
is of real type\cite{Ma:mec}. This is the approach which we consider in the
present paper. For the standard $SL_q(2)$ or Jones-polynomial solution $R$ of
the Quantum Yang-Baxter Equations (QYBE) we have the algebra
$BM_q(2)$\cite{Ma:exa}\cite{CWSSW:ten}
\ceqn{minkalg}{\vecu=\pmatrix{a&b\cr c&d},\quad ba=q^2ab,\quad
ca=q^{-2}ac,\quad d a=ad\\
bc=cb+(1-q^{-2})a(d-a),\quad d b=bd+(1-q^{-2})ab,\quad cd=d c+(1-q^{-2})ca}
which means that our braided matrix approach includes and extends the approach
of Carow-Watamura {\em et al.}  We define its $*$-algebra of co-ordinate
functions as $\R_q^{1,3}=BM_q(2)$. It has time direction central and the space
directions mutually non-commuting. Note that the braided approach replaces such
explicit relations by R-matrix formulae, which are easier to work with and more
general.

The chacteristic feature of this braided approach, as well as being always of a
general R-matrix form, is that the underlying objects
are not Hopf algebras or quantum groups in any usual sense.  Instead, we
require the new concept, introduced in \cite{Ma:rec}\cite{Ma:exa}\cite{Ma:bra},
of a {\em braided group}. This is like a quantum group but the coproduct map
$\Delta:B\to B\und\tens B$, say, is a homomorphism with respect to the {\em
braided tensor product} algebra in which the two tensor factors in $B\und\tens
B$ do not commute. Instead they enjoy mutual {\em braid statistics}, given in
our examples by $R$-matrices (and any $q$ or other parameters in them). Long
introductions to the general theory of braided groups are in \cite{Ma:introp}
and \cite{Ma:introm}. See also chapter~10 of my textbook on quantum
groups\cite{Ma:book}.

We begin in the preliminary Section~2 by recalling from
\cite{Ma:exa}\cite{Mey:new} the multiplicative and additive braided group
structures in the braided matrices $B(R)$. This indeed justifies this name for
the algebra (\ref{bramat}), for it corresponds when $q=1$ to the multiplication
and addition of usual matrices. For $\mink$ in the twistor viewpoint, the first
is needed to pick out a natural $q$-determinant or square-radius function
(which determines the quantum metric) and the second (which is due to U.
Meyer\cite{Mey:new}) gives the linear structure of spacetime (such as the
addition of 4-momentum).

The general framework of {\em braided coaddition} for linear braided group
structures was introduced in \cite{Ma:poi}, where we showed, quite generally,
how to build from this a Poincar\'e quantum group by a {\em bosonisation}
construction\cite{Ma:bos}\cite{Ma:poi}. We adjoin the $q$-Lorentz generators
by a kind of semidirect product which is different, however, from the
bicrossproducts needed for the $\kappa$-deformation mentioned above. The
abstract picture is summarized in the Appendix~A. The calculation for the
specific case of braided matrices $B(R)$ was done in \cite{MaMey:bra} and gave
the $q$-Poincar\'e (function algebra) quantum group in $R$-matrix form, for the
first time. We recall this is Section~3. We also explained at the end of
\cite{MaMey:bra} that once in $R$-matrix form, it is an easy matter to dualise
and obtain the Poincar\'e quantum group in enveloping algebra form. The modest
contribution in the present paper is to give the resulting formulae explicitly
for those who do not want to carry out the exercise. The result originates,
however, in \cite{MaMey:bra} and the formulae are very similar. This is covered
in Section~4. The standard $R$-matrix gives the Minkowski space $q$-Poincar\'e
quantum enveloping algebra which can be compared with explicit generators and
relations in \cite{OSWZ:def}.

Next we turn our attention to $q$-Euclidean space. There is a parallel
$R$-matrix theory for this, introduced (in the $R$-matrix form)  in
\cite{Ma:euc} based on the algebra $\bar A(R)$
\eqn{eucmat}{\vecx=\{x^i{}_j\},\quad  \eucmat }
in place of (\ref{bramat}). The standard $SU_q(2)$ or Jones polynomial
$R$-matrix gives an algebra $\bar M_q(2)$
\ceqn{eucalg}{\vecx=\pmatrix{a&b\cr c&d};\qquad ba=qab,\quad ca=q^{-1}ac,\quad
da=ad\\ bc=cb+(q-q^{-1})ad\qquad db=q^{-1}bd\quad dc=qcd  }
which is (by an accident) isomorphic to the usual FRT bialgebra $A(R)=M_q(2)$
of quantum matrices. The latter was proposed as $q$-Euclidean space in
\cite{CWSSW:ten} so, once again, the braided approach includes and extends that
 pioneering work in a general $R$-matrix form (\ref{eucmat}) introduced in this
context in \cite{Ma:euc}. We define the non-commuting co-ordinate functions as
$\R_q^4=\bar M_q(2)$. There is a natural unitary-like $*$-structure with
corresponding Euclidean-type norm defined by the natural $q$-determinant. Also
in \cite{Ma:euc} is a theory of twisting or `quantum Wick rotation' which says
that the $\bar A(R)$ algebra (\ref{eucmat}) is `gauge equivalent' (in a sense
generalising ideas of Drinfeld\cite{Dri:qua}) to our first $B(R)$ algebra
(\ref{bramat}). If we did not have a different $*$ then we would have something
strictly equivalent (via twisting) to the already-established $q$-Minkowski
space above, which would not be very interesting. Instead, we put the algebra
to good use as a complement to the $q$-Minkowski above\cite{Ma:euc}.

The treatment of (\ref{eucmat}) as an additive braided group (or braided
covector space) is in Section~5. We also give its multiplicative structure,
which is like the braided group $B(R)$ but the non-commutativity required does
not obey the Artin braid relations (or QYBE), so it is something a bit beyond
even the concept of a braided group. In Section~6 we recall the corresponding
q-Poincar\'e quantum group in R-matrix form as obtained by the general
bosonisation construction, applied now to (\ref{eucmat}). This is also
from\cite{Ma:euc}. Once again, we emphasized there the function algebra theory
and left the dualisation for the enveloping algebra $q$-Poincar\'e quantum
group as an easy exercise. For completeness, we do this now in Section~7. The
standard $R$-matrix gives the Euclidean group of motions and can be compared
with the $n=4$ case of computations based on $SO_q(n)$-covariant quantum planes
 in \cite{Fio:euc} and elsewhere. We also give explicitly the duality pairing
and, accordingly, the action of the Poincar\'e quantum group on spacetime
co-ordinates $B(R)$ and $\bar A(R)$. Note that we will use the term
`Poincar\'e' quantum group to cover all dimensions and signatures since our
construction is quite uniform.

Let us stress that this paper for the most part collects results and formulae
already obtained in some form in \cite{Ma:poi}\cite{MaMey:bra}\cite{Ma:euc}
modulo conventions and elementary dualisation. Nevertheless, it is hoped that a
self-contained account now will be useful as an overview and introduction. It
complements an extensive 50-page introduction to our `braided geometry'
approach to q-Minkowski space in \cite{Ma:varen}. The extensive further
literature on this topic, including works by the
author\cite{Ma:eps}\cite{Ma:star} and U. Meyer\cite{Mey:wav}, as well as by the
Munich and Berkeley groups are covered there. It could be said that the
underlying geometry is fairly well understood by now, though not the full story
regarding the $*$-structure and the construction of actual q-deformed quantum
field theories on these spacetimes, which remain a goal and motivation for the
subject. Such a q-deformation would surely have an application either as a tool
for regularising infinities (as poles at $q=1$)\cite{Ma:reg} or as a model of
quantum or other effects on the structure of geometry at the Planck
scale\cite{Ma:pla}\cite{Kem:fie}.

Finally, Section~8 contains some newer material. We point out that the
bosonisation point of view does suggest a natural solution to the problem of
what should be the correct axioms for the $*$-structure on our Poincar\'e
quantum groups. This makes contact with some preliminary ideas in the preprint
of Fiore\cite{Fio:man}. We also observe that one can avoid the dilatonic
extension needed in the constructions above by only partially bosonising. In
this case the Poincar\'e algebra is a braided group not a quantum group, but of
a fairly mild kind where the braid-statistics are given by a phase factor as in
\cite{Ma:csta}.

To be clear about normalisation, what we call the standard $SU_q(2)$ or
Jones-polynomial R-matrix is
\eqn{Rsl2}{R=\pmatrix{q&0&0&0\cr 0&1&q-q^{-1}&0\cr 0&0&1&0\cr 0&0&0&q}}
where the rows and columns label the two copies of $M_n$ in $M_n\tens M_n$,
wherein $R$ lives. Our formulae are quite general and not at all
limited to this particular R-matrix.

\section{Preliminaries I: $q$-Minkowski spaces in R-matrix form}

The braided-matrix approach to $q$-Minkowski space was presented at the
Zdikov Winter School in January 1993, at the Guadeloupe Spring
School\cite{Ma:carib} in May 1993 and in Clausthal\cite{Ma:clau} in
July 1993. So this preliminary section and the beginning of the next
will recall some of that basic material from
\cite{Ma:mec}\cite{Mey:new}.

The idea is that in classical geometry one can take the space of hermitian
matrices with norm given by the determinant, as Minkowski space. So let us
build a convincing $q$-deformation of the concept of a hermitian matrix. To be
a matrix, we need to be able to linearly add, and matrix multiply (with an
identity for the multiplication). In our
algebraic language we indeed have this on the braided matrices  $B(R)$ as a
{\em braided matrix comultiplication} and {\em counit}
\eqn{bramult}{\Delta\vecu=\vecu\tens\vecu,\quad \eps\vecu=\id}
where $\Delta$ is an algebra homomorphism $B(R)\to B(R)\und\tens B(R)$. The
{\em multiplicative braid statistics} needed for it to extend as an algebra
homomorphism are\cite{Ma:exa}
\eqn{mulstat}{\vecu''=\vecu\vecu';\qquad \mulstat}
whereby $\vecu''$ obeys the relations (\ref{bramat}) of $B(R)$ if
$\vecu,\vecu'$ do.

We also have a hermitian $*$-structure $u^i{}_j{}^*=u^j{}_i$ whenever $R$ is of
real-type\cite{Ma:mec}, obeying the axioms (introduced there) of a $*$-braided
group of real type,
\eqn{delta*}{\Delta\circ *=(*\tens *)\circ\tau\circ\Delta,\quad
\eps\circ*=\overline{\eps(\ )}.}
where $\tau$ is transposition. When there is an `inverse' or braided antipode
$S$ (which is not the case for braided matrix multiplication in $B(R)$ of
course) we demand also $*\circ S=S\circ *$.
The astute reader may wonder, by the way, about these axioms (\ref{delta*}).
They are quite different from the usual axioms of a Hopf $*$-algebra. Indeed,
usual hermitian matrices do not form a group under composition, which is why in
algebraic terms $\Delta$ is not a $*$-algebra homomorphism. But if $A,B$ are
hermitian then $(AB)^*=BA$, which is why $\Delta$ commutes with $*$ when we put
in the extra transposition $\tau$. This means that the algebraic notion of
$*$-braided groups or $*$-quantum groups with the transposed axiom
(\ref{delta*}) is useful even for $R=1$, where it allows us to view the space
of hermitian matrices (e.g. the mass-shell in Minkowski space) as a `group' in
this generalised sense.

Finally, we have, at least when $R$ is $q$-Hecke, a {\em braided
coaddition}\cite{Mey:new}
\eqn{braadd}{\Delta\vecu=\vecu\tens 1+1\tens\vecu,\quad\eps\vecu=0,\quad
S\vecu=-\vecu}
again extended as an algebra homomorphism $\Delta:B(R)\to B(R)\und\tens B(R)$.
This time we use Meyer's {\em additive braid statistics}\cite{Mey:new}
\eqn{addstat}{\vecu''=\vecu+\vecu';\qquad \addstat.}
When $R$ is of real type, we have again a $*$-braided group with our $*$ as
above and obeying the same axioms (\ref{delta*}). In this case the
transposition $\tau$ in the axioms would not be visible when $R=1$ because in
this case the coaddition $\Delta$ would be cocommutative. This is correct
because there is no problem adding hermitian matrices.

These basic properties amply justify the term `braided matrices' for the
algebra (\ref{bramat}). We also frequently write all the four $R$-matrices in
these equations on one side as `big' multiindex matrices ${\bf R}$, ${\bf R'}$
etc. where we consider $u^i{}_j=u_I$ as a vector with multiindex $I=(i_0,i_1)$.
This is how in \cite{Ma:exa} we wrote equations such as (\ref{bramat}), which
occur in other contexts too\cite{FRT:lie}, in `Zamalodchikov' or {\em braided
covector}  form
\eqn{bracovec}{\vecu_1\vecu_2=\vecu_2\vecu_1{\bf R}',\qquad
\vecu''=\vecu+\vecu';\qquad \vecu'_1\vecu_2=\vecu_2\vecu'_1{\bf R}}
for suitable ${\bf R}',{\bf R}$. They were introduced in
\cite{Ma:exa},\cite{Mey:new} under the names $\Psi',R_L$ respectively. Likewise
for (\ref{mulstat}) introduced in \cite{Ma:exa} as $\Psi$. We assume that the
reader can freely transfer back and forth between the multiindex ${\bf R},{\bf
R}'$ braided covector form and the previous matrix form (\ref{bramat}) etc.,
without imagining that they are anything but different notations. The reason we
use both forms is that the matrix form emphasizes the braided comultiplication
structure while the braided covector form emphasises the linear structure as
like a quantum plane. There is a detailed theory of braided covectors
$\Vhaj({\bf R}',{\bf R})$ introduced in \cite{Ma:poi}\cite{Ma:fre} which we can
now use for $B(R)$ as a linear braided space. This includes covariance and
Poincar\'e quantum groups\cite{Ma:poi}, a theory of braided-differentiation,
integration\cite{KemMa:alg}, epsilon tensor and electromagnetism\cite{Ma:eps},
and so on -- all the constructions we are familiar with for $\R^n$.

For our standard example where $R$ is (\ref{Rsl2}), we have the multiplicative
braid statistics\cite{Ma:exa}
\ceqn{minkmulstat}{a'  a=a a'+(1-q^2)b c',\quad a'  b=b  a',\quad  a'
c=ca'+(1-q^2)(d-a)  c'\\
a'  d=d  a'+(1-q^{-2})b  c',\quad  b'  a=a  b'+(1-q^2)b  (d'- a'),\quad
b'b=q^2b  b',\quad{\rm etc.} }
for $\pmatrix{a''&b''\cr c''&d''}=\pmatrix{a&b\cr c&d}\pmatrix{a'&b'\cr c'&d'}$
to obey the relations (\ref{minkalg}), and additive braid
statistics\cite{Mey:new}
\ceqn{minkaddstat}{a'  a=q^2a a',\quad a'  b=b  a',\quad a'
c=q^2ca'+(q^2-1)ac',\quad c'a=ac'\\
a'  d=d  a'+(q^2-1)b  c' +q^{-2}(q^2-1)^2aa',\quad  b'  a=q^2a
b'+(q^2-1)ba',\quad b'b=q^2b  b',\quad {\rm etc.} }
for $\pmatrix{a''&b''\cr c''&d''}=\pmatrix{a&b\cr c&d}+\pmatrix{a'&b'\cr
c'&d'}$ to obey (\ref{minkalg}).

The braided comultiplication is needed (as for usual matrices) to fix the {\em
braided determinant} as group-like or `multiplicative'. It comes out as
\cite{Ma:exa}
\eqn{bradet}{\und\det(\vecu)=ad-q^2cb}
and is central, as well as bosonic with respect to the multiplicative braid
statistics. We use is as a square-distance function on $BM_q(2)$.

Finally, we have the hermitian $*$,
\eqn{herm*}{\pmatrix{a*&b*\cr c*&d*}=\pmatrix{a&c\cr b&d}}
so that these matrices are indeed naturally hermitian, and hence obviously
provide a natural definition $\mink=BM_q(2)$, as explained in \cite{Ma:mec}.
The $*$-structure is needed to determine what should be the `real' or
self-adjoint space-time co-ordinates under $*$. The required linear
combinations are
\eqn{xyzt}{t={q^{-1}a+qd\over 2},\quad x={b+c\over 2},\quad y={b-c\over
2i},\quad z={d-a\over 2}}
and the braided determinant becomes\cite{Ma:mec}
\eqn{detxyzt}{\und{\rm det}(\vecu)={4q^2\over(q^2+1)^2}t^2-q^2x^2-q^2 y^2-
{2(q^4+1)q^2\over (q^2+1)^2}z^2+\left({q^2-1\over q^2+1}\right)^2{2q}
tz}
which justifies indeed the interpretation as Minkowski length from the braided
approach. From it we can extract a
{\em quantum metric} tensor by braided differentiation\cite{Ma:eps} which in
our matrix basis is
\eqn{minkmet}{ \eta^{IJ}=\pmatrix{q^{-2}-1&0&0&1\cr 0&0&-q^{-2}&0\cr
0&-1&0&0\cr 1&0&0&0}; \quad \und\det(\vecu)=(1+q^{-2})^{-1}\eta^{IJ}u_Ju_I.}

\section{Minkowski $q$-Poincar\'e quantum group in function algebra form}

One of the first consequences of writing the braided matrices $B(R)$ algebra
(\ref{bramat}) in the braided covector form (\ref{bracovec}) is that we know at
once from \cite{Ma:poi} what its associated Poincar\'e quantum group looks
like. For in \cite{Ma:poi} was introduced a completely general $R$-matrix
construction for such objects based on a theory of `bosonisation' in
\cite{Ma:bos}. The formulae are as follows. For the `rotation group' we take
the usual FRT bialgebra\cite{FRT:lie} $A({\bf R})$ with generator
$\Lambda^I{}_J$ say, and relations
\eqn{FRTL}{ {\bf R}\Lambda_1\Lambda_2=\Lambda_2\Lambda_1 {\bf R},\quad \Delta
\Lambda=\Lambda\tens\Lambda,\quad \eps\Lambda=\id}
to which we add relations needed to give us a Hopf algebra with antipode. They
include such things as a metric relation
\eqn{lambdaeta}{\Lambda_1\Lambda_2\eta_{21}=\eta_{21},\quad {\rm i.e.},\quad
S\Lambda^I{}_J=\Lambda^A{}_B\eta^{BI}\eta_{AJ}}
where we let $\eta_{IJ}$ be the transposed inverse of $\eta^{IJ}$. This is the
vectorial approach to the Lorentz quantum group in \cite{Mey:new}. Next, we
extend this by adjoining a central group-like element $\dila$ (the dilaton) in
such a way that our braided covectors are fully covariant under the extended
transformation
\eqn{veccov}{ \vecu\to\vecu\Lambda\dila}
with additive braid statistics correctly induced by this covariance. We need
$\dila$ to achieve this because ${\bf R}$ as given is not in the `quantum group
normalisation' and we have to compensate for this\cite{Ma:poi}. We now use the
braided covector generators $\vecu$ for the momentum sector of our
$q$-Poincar\'e quantum group, and denote them as such by $\vecp$ to avoid
confusion with spacetime co-ordinates. We now make a coalgebra semidirect
coproduct by the above coaction. At the same time, the dual quasitriangular
structure\cite{Mey:new}  (universal R-matrix functional) of our extended
Lorentz quantum group converts this coaction to an action, and we make an
algebra semidirect product by this. The theory of bosonisation ensures that the
result is necessarily a Hopf algebra. Some of this general theory is recalled
in Appendix~A. All we need to know for now is the resulting $R$-matrix
formula\cite{Ma:poi}
\ceqn{vecpoi}{ \vecp_1\vecp_2=\vecp_2\vecp_1{\bf R}',\quad
\vecp_1\Lambda_2=\lambda\Lambda_2\vecp_1{\bf R},\quad \vecp \dila
=\lambda^{-1}\dila \vecp,\quad [\Lambda,\dila ]=0\\
 \Delta \dila =\dila \tens \dila,\quad \eps\dila=1,\quad S\dila=\dila^{-1}\\
 \Delta \vecp=\vecp\tens \Lambda \dila +1\tens\vecp,\quad \eps\vecp=0,\quad
S\vecp=-\vecp \dila^{-1}\Lambda^{-1}}
where $\lambda$ is the {\em quantum group normalisation constant} of ${\bf R}$
defined such that $\lambda{\bf R}$ is in the quantum group
normalisation\cite{Ma:lin}. Some special cases of this construction for
$ISO_q(n)$, etc., were first studied in \cite{SWW:inh}, though without the
above general construction. Finally, no proposal for a $q$-Poincar\'e quantum
group is complete without a covariant action of it on the $q$-spacetime
co-ordinates. The general construction\cite{Ma:poi} introduced just this, in
the coaction form
\eqn{veccoact}{ \vecu\to \vecu\Lambda\dila+\vecp}
extended as an algebra homomorphism. The $\vecu$ commute with the Poincar\'e
generators. This covariance was one of the main achievements of the braided
approach in \cite{Ma:poi} and is a general feature of the bosonisation theory
recalled in Appendix~A.

This is the {\em vectorial form} of the $q$-Poincar\'e quantum group when we
apply it to (\ref{bramat}) in the form (\ref{bracovec}). There is also a {\em
spinorial form} obtained when we unwind the above construction back in terms of
$R$ rather than ${\bf R}',{\bf R}$. Firstly, we replace the quantum group
$A({\bf R})$ by $A(R)$ and make this into a Hopf algebra $A$. Two copies of it
are needed to play the role of $A({\bf R})$ above, with generators
$\vecs,\vect$ say, forming a double cross product Hopf algebra $A\dcross
A$\cite{Ma:mor}
\ceqn{spinlor}{R\vecs_1\vecs_2=\vecs_2\vecs_1R,\quad\FRT,\quad
R\vect_1\vecs_2=\vecs_2\vect_1 R\\
\Delta\vecs=\vecs\tens\vecs,\quad \Delta\vect=\vect\tens\vect,\quad
\eps\vecs=\id,\quad \eps\vect=\id}
among the further relations needed to give $\vecs,\vect$ antipodes. The
vectorial form is realised in terms of the spinorial form by
\eqn{lambdast}{ \Lambda^I{}_J=(S s^{j_0}{}_{i_0}) t^{i_1}{}_{j_1}}
and the covariance (\ref{veccov}) takes the form cf\cite{CWSSW:lor}
\ceqn{minkcov}{ \vecu\to \vecs^{-1}\vecu\vect\dila.}
We again adjoin a central grouplike dilaton $\dila$ to take care of the fact
that $R$ is not in the quantum group normalisation. The general construction of
the quantum group $A\dcross A$ was developed in \cite[Sec. 4]{Ma:poi}, as well
as the isomorphism with Drinfeld's quantum double\cite{Dri} and the
identification (relevant later) as a twisting of $A\tens A$. It connects in the
$SU_q(2)\dcross SU_q(2)$ case with the proposal for $q$-Lorentz group in
\cite{CWSSW:lor} and with the proposal as the quantum double in
\cite{PodWor:def}. There is a $*$-structure
\eqn{spin*}{ s^i{}_j{}^*=St^j{}_i,\quad t^i{}_j{}^*=Ss^j{}_i,\quad
\dila^*=\dila}
as in \cite{CWSSW:lor}, which now works for general $R$ of
real-type\cite{Ma:poi}. The diagonal case $\vect=\vecs$ (without the dilaton)
is covariance under the spacetime spinor rotation group $SU_q(2)$ and preserves
the multiplicative structure, the distance function $\und\det(\vecu)$ and the
time co-ordinate $t$\cite{Ma:exa}.

The formulae (\ref{vecpoi}) then become the spinorial $q$-Poincar\'e quantum
group\cite{MaMey:bra}
\ceqn{spinpoi}{R_{21}\vecp_1R\vecp_2=\vecp_2R_{21}\vecp_1R,\quad
\vecp_1\vecs_2=\vecs_2 R^{-1}\vecp_1 R,\quad
\vecp_1\vect_2=\lambda\vect_2R_{21}\vecp_1 R\\
\vecp \dila =\lambda^{-1}\dila \vecp,\quad [\vecs,\dila]=[\vect,\dila]=0,\quad
\Delta\dila=\dila\tens\dila,\quad \eps\dila=1,\quad S\dila=\dila^{-1}\\
\Delta \vecp=\vecp\tens\vecs^{-1}(\ )\vect \dila +1\tens\vecp,\quad
\eps\vecp=0,\quad S\vecp=-\vecp S(\vecs^{-1}(\ )\vect\dila)}
where $\vecs^{-1}(\ )\vect$ has a space for the matrix indices of $\vecp$ to be
inserted. The constant $\lambda$ is the square of the quantum group
normalisation constant of $R$. Its value for our standard example is
$\lambda=q^{-1}$. One can again derive this construction by the more abstract
bosonisation construction (\ref{cobos}) in Appendix~A, knowing only the full
covariance (\ref{minkcov}). The two methods give the same answer. Finally, the
coaction of the $q$-Poincar\'e quantum group on the spacetime co-ordinates is
the algebra homomorphism
\eqn{spincoact}{\vecu\to \vecs^{-1}\vecu\vect\dila+\vecp}
where the $\vecu$ commute with the Poincar\'e generators.

Recall that we also have a multiplicative braided group structure on $B(R)$.
The `mass shell' in $q$-Minkowski space is the $*$-braided group $BSU_q(2)$ in
\cite{Ma:exa} with (braided) antipode. If we bosonise by the above rotational
$SU_q(2)$ covariance, we obtain this time\cite{Ma:mec} the quantum double of
$SU_q(2)$ which, as noted above, is isomorphic to the $q$-Lorentz group
$SU_q(2)\dcross SU_q(2)$ in spinorial form. Such an isomorphism has no
classical counterpart, being singular at $q=1$.

\section{Minkowski $q$-Poincar\'e quantum group in enveloping algebra form}

Once we have our Poincar\'e quantum groups in $R$-matrix form, it is a pretty
easy exercise to dualise them and find the corresponding Poincar\'e quantum
group enveloping algebras. We just dualise the individual pieces of the
bosonisation or semidirect construction. We suppose that the dual quantum group
to the Lorentz quantum group function algebra can be put in FRT form with
generators ${\bf L}^\pm$ say and relations\cite{FRT:lie}
\eqn{Lpm}{\vecL^\pm_1\vecL^\pm_2{\bf R}={\bf R}\vecL^\pm_2\vecL^\pm_1,\quad
\vecL^-_1\vecL^+_2{\bf R}={\bf R}\vecL^+_2\vecL^-_1,\quad \Delta
\vecL^\pm=\vecL^\pm\tens \vecL^\pm,\quad \eps \vecL^\pm=\id}
among others needed in the dual. We also define the dual of the dilaton
generator $\dila$ to be the generator $\xi$ of the enveloping algebra $U_q(1)$,
which we define with a non-standard universal $R$-matrix\cite{LyuMa:fou}. We
extend our Lorentz enveloping algebra by this. Finally, for the dual of our
previous momentum co-ordinates $p_I$ as a braided covector space $\Vhaj({\bf
R}',{\bf R})$ we need look no further than the corresponding {\em braided
vector algebra} $V({\bf R}',{\bf R})$ which is arranged carefully in
\cite{Ma:poi}\cite{Ma:fre} to be the braided group dual to the braided
covectors. Its generators $P^I$, say, have upper indices. The pairing between
our objects in the last section and our new enveloping algebra generators is
then
\eqn{vecpair}{ \<\Lambda_1,\vecL^+_2\>=\lambda{\bf R},\quad
\<\Lambda_1,\vecL^-_2\>=\lambda^{-1}{\bf R}_{21}^{-1},\quad
\<P^I,p_J\>=\delta^I{}_J,\quad \<\dila,\xi\>=1}
with the trivial pairing provided by the counits between the different quantum
or braided groups.

Using this pairing we can then deduce the Poincar\'e quantum group in
enveloping algebra form in the standard way by dualising. The resulting
structure is\cite{Ma:poi}
\ceqn{vecpoienv}{ \vecP_1\vecP_2={\bf R}'\vecP_2\vecP_1,\quad \vecL_1^+
\vecP_2=\lambda^{-1}{\bf R}_{21}^{-1}\vecP_2\vecL^+_1,\quad
\vecL_1^-\vecP_2=\lambda {\bf R}\vecP_2\vecL^-_1,\quad
\lambda^\xi\vecP=\lambda^{-1}\vecP\lambda^{\xi}\\
{}[\xi,\vecL^\pm]=0,\quad \Delta\xi=\xi\tens 1+1\tens\xi,\quad \eps\xi=0,\quad
S\xi=-\xi\\
\Delta\vecP=\vecP\tens 1+\lambda^\xi\vecL^-\tens\vecP,\quad \eps \vecP=0,\quad
S\vecP=-\lambda^{-\xi} (S\vecL^-)\vecP.}
We can also obtain the above formulae by starting with covariance under the
Hopf algebra obtained from (\ref{Lpm}) and dilatonic extension, namely
\eqn{Lpmcov}{\vecL^+_1\la\vecP_2=\lambda^{-1}{\bf R}^{-1}_{21}\vecP_2,\quad
\vecL^-_1\la\vecP_2=\lambda{\bf R}\vecP_2,\quad
\lambda^\xi\la\vecP=\lambda^{-1}\vecP}
and applying it to the enveloping-type bosonisation theorem (\ref{bos}) in
Appendix~A. Finally, we can dualise our coaction (\ref{veccoact}) on the
spacetime co-ordinates by evaluating the relevant part of the output of the
coaction against our $q$-Poincar\'e enveloping algebra generators. It obviously
becomes now an action of the $q$-Poincar\'e enveloping algebra
\eqn{vecact}{\vecL^+_1\la\vecu_2=\vecu_2\lambda{\bf R}_{21},\quad
\vecL^-_1\la\vecu_2=\vecu_2\lambda^{-1}{\bf
R}^{-1},\quad\lambda^\xi\la\vecu=\lambda\vecu,\quad P^I\la u_J=\delta^I{}_J}
necessarily making our position co-ordinates $\vecu$ a module algebra. This
means that it extends to products using the coproduct in (\ref{vecpoi}). For
example, we see easily that
\eqn{Pact}{P^I\la (u_{I_1}\cdots u_{I_m})=u_{J_2}\cdots u_{J_m}\left[m;{\bf
R}^{-1}_{21}\right]^{IJ_2\cdots J_m}_{I_1I_2\cdots I_m}}
where $[m;{\bf R}^{-1}_{21}]$ is a {\em braided integer matrix} as introduced
in \cite{Ma:fre} in the theory of braided differentiation on linear braided
spaces. It means that the realisation of the $P^I$ generators acting on our
braided spacetime is exactly by braided differentiation. We concentrated in
\cite{Ma:fre} on braided differentials $\del^I$ with matrix $[m;{\bf R}]$; the
other case $\bar\del^I$ say, is equivalent by a symmetry principle\cite[Sec.
V]{Ma:fre}. If $\Vhaj({\bf R}',{\bf R})$ is a braided covector algebra then so
is $\Vhaj({\bf R}',{\bf R}^{-1}_{21})$. It lives in the `opposite' category
with inverse-transposed braiding. So $P^I=\bar\del^I$ is our natural
realisation of the momentum sector.

One of the main results in \cite{Ma:fre} is that $\del^I$ (and hence also
$\bar\del^I$) always obey the relations of the braided vector algebra $V({\bf
R}',{\bf R})$, which confirms their role now as representing the $P^I$. The
$\del^I$ extend to products with the inverse braiding $\Psi^{-1}$ and the
$\bar\del^I$ with the usual braiding $\Psi$, i.e. by a {\em braided Leibniz
rule}\cite{Ma:fre}. If we consider  $u_I$ also as an operator by left
multiplication then the corresponding braided-Leibniz rules are expressed in
general as the commutation relations
\eqn{leib}{ \del_1\vecu_2-\vecu_2{\bf R}_{21}\del_1=\id,\quad
\bar\del_1\vecu_2-\vecu_2{\bf R}^{-1}\bar\del_1=\id.}
This generalised and included considerations for specific
algebras\cite{PusWor:twi}\cite{Kem:sym}\cite{OSWZ:def} to the completely
general setting of differentiation on any braided covector space.

This is the vectorial form of the $q$-Poincar\'e enveloping algebra according
to the general construction for any braided covector space \cite{Ma:poi}. Next
we note that when there is quantum metric $\eta$, it provides\cite{Mey:new} an
isomorphism between the braided vectors and covectors as braided groups. This
is in fact the abstract definition of a quantum metric in braided geometry.
Hence in this case we don't have to work with braided covectors (though that is
more canonical) but can instead use the quantum metric to lower the indices and
work everywhere with braided covectors. To do this one needs the standard
quantum metric identities
\eqn{metricid}{ \eta_{IA}{\bf R}^{-1}{}^A{}_J{}^K{}_L=\lambda^2 {\bf
R}^A{}_I{}^K{}_L\eta_{AJ},\quad \eta_{KA}{\bf R}^I{}_J{}^A{}_L=\lambda^{-2}{\bf
R}^{-1}{}^I{}_J{}^A{}_K\eta_{AL},\quad {\rm etc.}}
deduced automatically from (\ref{lambdaeta}) by evaluating it against the
universal R-matrix functional. If we make this change of basis to new covector
generators $P_I=\eta_{IA}P^A$, then the vectorial $q$-Poincar\'e enveloping
algebra with braided covectors $\vecP$ is clearly
\ceqn{lowvecpoienv}{ \vecP_1\vecP_2=\vecP_2\vecP_1{\bf R}',\quad \vecL_1^+
\vecP_2=\lambda\vecP_2{\bf R}_{21}\vecL^+_1,\quad \vecL_1^-\vecP_2=\lambda^{-1}
\vecP_2{\bf R}^{-1}\vecL^-_1,\quad
\lambda^\xi\vecP=\lambda^{-1}\vecP\lambda^{\xi}\\
{}[\xi,\vecL^\pm]=0,\quad \Delta\xi=\xi\tens 1+1\tens\xi,\quad \eps\xi=0,\quad
S\xi=-\xi\\
\Delta P_I=P_I\tens 1+ \lambda^\xi SL^-{}^A{}_I\tens P_A,\quad \eps
\vecP=0,\quad SP_I=-\lambda^{-\xi} (S^2 L^-{}^A{}_I)P_A.}
The structure is a bosonisation according to (\ref{bos}) in Appendix~A by
$\vecL^\pm$ acting as in (\ref{vecact}) since the $P_I$ transform like the
$u_I$ under these, and $\lambda^\xi\la P_I=\lambda^{-1}P_I$. The pairing
(\ref{vecpair}) with the function algebra Poincar\'e generators from the last
section and the action (\ref{vecact}) on the spacetime co-ordinates of course
get modified to
\eqn{lowmom}{\<P_I,p_J\>=\eta_{IJ},\quad P_I\la u_J=\eta_{IJ}.}

One can do the same lowering of indices for the braided differentiation
operators to $\del_I=\eta_{IA}\del^A$ and $\bar\del_I=\eta_{IA}\bar\del^A$, in
which case these obey the braided covector relations and braided-Leibniz rule
\ceqn{lowdif}{ \del_1\del_2=\del_2\del_1{\bf R}',\quad
\bar\del_1\bar\del_2=\bar\del_2\bar\del_1{\bf R}'\\
\del_1\vecu_2-\lambda^{-2}\vecu_2\del_1{\bf R}^{-1}_{21}=\eta,\quad
\bar\del_1\vecu_2-\lambda^2\vecu_2\bar\del_1{\bf R}=\eta}
by the same elementary rearrangements using the quantum metric identities
(\ref{metricid}). So our lower-index momentum generators are represented on
spacetime by these lowered index braided differentials $P_I=\bar\del_I$.

This is the vectorial Poincar\'e group enveloping algebra with lower indices.
As in the start of Section~3, all the above is quite general. Now we apply it
to the specific form of ${\bf R}',{\bf R}$ in (\ref{bracovec}) in Section~2,
and unwind in terms of $R$ and $\vecu$ as a braided matrix $B(R)$. This gives
the $q$-Poincar\'e enveloping algebra in spinorial form. So we let $\vecl^\pm$
be dual to the generator $\vecs$ and $\vecm^\pm$ dual to the generator $\vect$
for each copy of a quantum enveloping algebra $H$ dual to $A$, with pairing
given by $R$ as in \cite{FRT:lie}. In the standard case they are each copies of
$U_q(su_2)$ in FRT form. The required Lorentz enveloping algebra is a {\em
double cross coproduct} $H\codcross H$ and does {\em not} contain the copies of
$H$ as sub-Hopf algebras. Instead they have a twisted
coproduct\cite{ResSem:mat}
\ceqn{spinlorenv}{\vecl^\pm_1\vecl^\pm_2R=R\vecl^\pm_2\vecl^\pm_1,\quad
\vecl^-_1\vecl^+_2R=R\vecl^+_2\vecl_1^-,\quad
\vecm^\pm_1\vecm^\pm_2R=R\vecm^\pm_2\vecm^\pm_1\\
\vecm^-_1\vecm^+_2R=R\vecm^+_2\vecm_1^-,\quad
[\vecl_1^\pm,\vecm_2^\pm]=[\vecl^\pm_1,\vecm^\mp_2]=0\\
\Delta\vecl^\pm=\CR^{-1}(\vecl^\pm\tens\vecl^\pm)\CR,\quad
\Delta\vecm^\pm=\CR^{-1}(\vecm^\pm\tens\vecm^\pm)\CR,\quad
\eps\vecl^\pm=\id,\quad \eps\vecm^\pm=\id}
where $\CR$ is the universal R-matrix or quasitriangular structure\cite{Dri} of
$H$ viewed as $\CR_{21}\in H\codcross H$. There is of course an appropriate
antipode and $*$, correspondingly twisted. We can realise the vectorial form
above in terms of these spinorial generators by cf. \cite{Mey:new}
\eqn{Llm}{
L^+{}^I{}_J=(\vecl^-\vecm^+)^{i_1}{}_{j_1}((S_0^{-1}\vecm^+)
(S_0^{-1}\vecl^+))^{j_0}{}_{i_0},\quad
L^-{}^I{}_J=(\vecl^-\vecm^-)^{i_1}{}_{j_1}
((S_0^{-1}\vecm^+)(S_0^{-1}\vecl^-))^{j_0}{}_{i_0}}
where $S_0$ is the usual `matrix inverse' antipode of $H$. The covariance of
the spacetime generators is expressed now as $B(R)$ a module algebra
under\cite{Ma:euc}
\ceqn{minkact}{\vecl^+_1\la \vecu_2=\lambda^{-\h}R^{-1}_{21}\vecu_2,\quad
\vecl^-_1\la\vecu_2=\lambda^\h R\vecu_2,\quad \vecm^+_1\la \vecu_2=\vecu_2
\lambda^\h R_{21},\quad \vecm^-_1\la\vecu_2=\vecu_2 \lambda^{-\h}R^{-1}}
which is also the action we use on the lowered momentum generators
$P_I=P^{i_0}{}_{i_1}$ regarded now as a braided matrix $B(R)$. We add the
dilaton $\xi$ as before. Semidirect product by this action (and the induced
semidirect coproduct) according to the bosonisation construction (\ref{bos}) in
Appendix~A, or just proceeding from (\ref{lowvecpoienv}) for the specific form
of ${\bf R}$ from (\ref{bracovec}), yields the Poincar\'e quantum group in
spinorial form as
\ceqn{spinpoienv}{R_{21}\vecP_1R\vecP_2=\vecP_2R_{21}\vecP_1R,\quad
\vecl^+_1\vecP_2\vecl^-_2=\lambda^{-\h}R^{-1}_{21}\vecP_2\vecl^-_2\vecl^+_1,
\quad \vecl^-_1\vecP_2\vecl^-_2=\lambda^\h R\vecP_2\vecl^-_2\vecl^-_1\\
 \vecm^+_1 \vecP_2\vecl^-_2=\vecP_2\vecl^-_2\lambda^\h R_{21}\vecm^+_1,
\quad\vecm^-_1\vecP_2\vecl^-_2=\vecP_2\vecl^-_2
\lambda^{-\h}R^{-1}\vecm^-_1,\quad
\lambda^\xi\vecP=\lambda^{-1}\vecP\lambda^\xi\\
{}[\xi,\vecl^\pm]=[\xi,\vecm^\pm]=0,\quad \Delta\xi=\xi\tens 1+1\tens\xi,
\quad\eps\xi=0,\quad S\xi=-\xi\\
\Delta \vecP=\vecP\tens 1+\lambda^\xi \vecl^-\vecm^+(\
)(S_0\vecm^-)(S_0\vecl^-)\tens\vecP,\quad \eps \vecP=0\\
 S\vecP=-\lambda^{-\xi}S_0\left(\vecm^+\vecl^-(\
)(S_0\vecl^-)(S_0\vecm^-)\right)\vecP}
where $(\ )$ is a space for the matrix entries of $\vecP$ to be inserted, and
$S_0$ is the usual matrix antipode in either copy of $H$. As far as I know,
these are new formulae presented here for the first time, although following
directly from \cite{Ma:poi} as above, or by dualisation of \cite{MaMey:bra}.
Finally, the action on the spacetime generators becomes, in this spinorial form
\ceqn{spinact}{\vecl^+_1\la \vecu_2=\lambda^{-\h}R^{-1}_{21}\vecu_2,\quad
\vecl^-_1\la\vecu_2=\lambda^\h R\vecu_2,\quad \vecm^+_1\la \vecu_2=\vecu_2
\lambda^\h R_{21},\quad \vecm^-_1\la\vecu_2=\vecu_2 \lambda^{-\h}R^{-1}\\
\lambda^\xi\la\vecu=\lambda\vecu,\quad \vecP_1\la \vecu_2=\eta.}
The action of the lowered $P_I$ is by the lowered $\bar\del_I$ automatically
obeying relations and braided-Leibniz rule
\eqn{lowminkdif}{R_{21}\bar\del_1 R\bar\del_2=\bar\del_2 R_{21}\bar\del_1
R,\quad R^{-1}\bar\del R\vecu_2-\lambda^2\vecu_2R_{21}\bar\del_1R=R^{-1}\eta\uo
R\eta\ut}
where $\eta=\eta\uo\tens\eta\ut$ is
$\eta_{IJ}=\eta^{i_0}{}_{i_1}{}^{j_0}{}_{j_1}$ as an element of $M_n\tens M_n$.
Some authors\cite{AzcRod:dif} have recently considered Leibniz rules for
$\mink$ of similar form on the left hand side. We stress, however, that this
is not a new equation but just the usual (\ref{leib}) or (\ref{lowdif}) for the
specific ${\bf R}$ introduced in \cite{Mey:new}. The rearrangement between the
two notations is just as in (\ref{addstat}) and (\ref{bracovec}), and quite
routine since\cite{Ma:skl}. Moreover, our treatment here works for general $R$
of $q$-Hecke type, with $\lambda$ the square of its quantum group normalisation
constant.

\section{Preliminaries II: $q$-Euclidean spaces in R-matrix form}

Now we make the same constructions as above for the $\bar A(R)$ algebra
(\ref{eucmat})  which was introduced in this context in\cite{Ma:euc}. As an
algebra it is exactly gauge equivalent to (\ref{bramat}) by the comodule
algebra twisting theory in \cite{Ma:euc}. We take, however, a non-hermitian
$*$-structure which ends up in the standard case as more like Euclidean space
than Minkowski.

To describe the multiplicative structure, we need to generalise the concept of
braided group $B$ slightly. We still require the coproduct $\Delta: B\to
B\und\tens B$ to be an algebra homomorphism, but don't insist that $\und\tens $
is a braided tensor product. The most general concept (which is certainly
general enough, but probably too general) is that $B\und\tens B$ should now be
some algebra which contains the two copies of $B$ as subalgebras, and uniquely
factorises into them in the sense that the map $B\tens B\to B\und\tens B$ given
by including the subalgebras and multiplying, is a linear isomorphism. This is
the algebra part of the theory of Hopf algebra factorisations in \cite[Sec.
3.2]{Ma:phy}, and more recently in \cite{Ma:prague} and elsewhere.

In this sense, $\bar A(R)$ does have a matrix comultiplication
\eqn{eucmult}{\Delta\vecx=\vecx\tens\vecx,\quad \eps\vecx=\id}
where $\Delta:\bar A(R)\to \bar A(R)\und\tens\bar A(R)$ is an algebra
homomorphism. So this is like a braided group (and indeed is gauge equivalent
to (\ref{bramult})) but the non-commutation relations describing $\und\tens$ do
not obey the Artin braid relations of QYBE; rather some more general (but not
completely general) algebra factorisation. The {\em multiplicative statistics}
are\cite{Ma:euc}
\eqn{bramulstat}{\vecx''=\vecx\vecx'; \qquad \vecx_1'\vecx_2=\vecx_2
R^{-1}\vecx_1'}
whereby $\vecx''$ obeys the same relations (\ref{eucmat}) if $\vecx,\vecx'$ do.

We also have a natural $*$-structure on $\bar A(R)$ whenever its more familiar
cousin $A(R)$ has a $*$-structure with real-type  universal R-matrix
functional. By a theorem in \cite{Ma:euc} we can take the same operation on the
$\vecx$ as we would on the $\vect$ of $A(R)$ in this case. It tends to be of
the type which is unitary in the quotient Hopf algebra when $R$ is of
real-type, namely of the form
\eqn{uni*}{x^i{}_j{}^*=\eps_{ai}x^a{}_b\eps^{bj}}
where $\eps^{ij}$ is invariant under the quotient Hopf algebra and $\eps_{ij}$
the transposed inverse.

Finally, we have a normal braided coaddition, under which $\bar A(R)$ remains
an additive braided group, at least when $R$ is $q$-Hecke\cite{Ma:euc}. We have
\eqn{eucadd}{\Delta\vecx=\vecx\tens 1+1\tens\vecx,\quad\eps\vecx=0,\quad
S\vecx=-\vecx}
extended as a braided group with {\em additive braid statistics}\cite{Ma:euc}
\eqn{baraddstat}{\vecx''=\vecx+\vecx';\qquad
\vecx'_1\vecx_2=R\vecx_2\vecx_1'R.}
This coaddition is typically compatible with the $*$ to give a $*$-braided
group as in (\ref{delta*}).

Clearly, we can also move the $R$'s to one side and write $\bar A(R)$ as a
braided covector algebra $x^{i_0}{}_{i_1}=x_I$ with (now in covector form as in
(\ref{bracovec})) the relations
\eqn{eucbracovec}{\vecx_1\vecx_2=\vecx_2\vecx_1{\bf R}',\qquad
\vecx''=\vecx+\vecx';\quad \vecx'_1\vecx_2=\vecx_2\vecx'_1{\bf R}}
for suitable ${\bf R}',{\bf R}$. These are given explicitly in \cite{Ma:euc}.
Thus we can equally well write this algebra with its linear braided group
structure in the form $\Vhaj({\bf R}',{\bf R})$ needed for our general
constructions in \cite{Ma:poi}\cite{Ma:fre}.

For our standard example where $R$ is (\ref{Rsl2}), we have the multiplicative
statistics\cite{Ma:euc}
\ceqn{eucmulstat}{a'a=q^{-1}aa'+(q^{-1}-q)bc',\quad a'b=ba',\quad
a'c=q^{-1}ca'+(q^{-1}-q)dc',\quad a'd=da'\\
b'a=q^{-1}ab'+(q^{-1}-q)bd',\quad b'b=bb',\quad
b'c=q^{-1}cb'+(q^{-1}-q)dd',\quad {\rm etc.}}
whereby the matrix product of the unprimed and primed matrices obeys the same
relations (\ref{eucalg}). We also have the additive braid statistics
\ceqn{eucaddstat}{a'a=q^2 aa',\quad a'b=qba',\quad a'c=qca'+(q^2-1)ac',\quad
a'd=da'+(q-q^{-1})bc'\\
b'a=qab'+(q^2-1)ba',\quad b'b=q^2bb',\quad
b'c=cb'+(q-q^{-1})(da'+ad')+bc'(q-q^{-1})^2, \quad {\rm etc.}}
\note{c'c=q^2 cc',\quad d'd=q^2 dd', \quad b'd=qdb'+(q^2-1)bd',\quad
c'a=qac',\quad c'b=bc',\quad c'd=qdc',\quad
d'b=qbd',\quad d'c=qcd'+(q^2-1)dc',\quad  d'a=ad'+(q-q^{-1})bc'}
whereby the sum of the unprimed and primed matrices obey the same relations
(\ref{eucalg}).

{}From the generalised comultiplication, or from the close connection of $\bar
M_q(2)$ with usual quantum matrices, one has a natural quantum determinant
\eqn{eucdet}{\und{\det}(\vecx)=ad-qcb}
which is central, as well as bosonic with respect to the multiplicative. We use
is as the square-distance function on $\bar M_q(2)$.

Finally, it is easy to see that the standard $2\times 2$ quantum matrices
$M_q(2)$ for real $q$ are a $*$-bialgebra with
\eqn{eucl*}{\pmatrix{a^*&b^*\cr c^*&d^*}=\pmatrix{d&-q^{-1}c\cr -qb&a}}
and hence from theory in \cite{Ma:euc} we know that our algebra $\bar M_q(2)$
is also a $*$-algebra with the operation (\ref{eucl*}). Under the coaddition it
forms a $*$-braided group obeying (\ref{delta*}).  Such a `unitary' form
provides us a natural definition $\R_q^4=\bar M_q(2)$, as explained in
\cite{Ma:euc}. The $*$-structure determines `real' or self-adjoint (under $*$)
spacetime co-ordinates
\eqn{euctxyz}{t={a-d\over 2\imath},\quad x={c-qb\over 2}, \quad y={c+qb\over
2\imath},\quad z={a+d\over 2}}
and the $q$-determinant above becomes
\eqn{euclmetric}{ \und{\det}(\vecx)=({1+q^{2}\over
2})t^2+x^2+y^2+({1+q^{2}\over 2})z^2}
which justifies the interpretation as Euclidean length in this
approach\cite{Ma:euc}. From it we can extract a
{\em quantum metric} tensor by braided differentiation\cite{Ma:eps} which in
our matrix basis is
\eqn{eucmet}{ \eta^{IJ}=\pmatrix{0&0&0&1\cr 0&0&-q^{-2}&0\cr 0&-1&0&0\cr
1&0&0&0}; \quad \und\det(\vecx)=(1+q^{-2})^{-1}\eta^{IJ}x_Jx_I.}

\section{Euclidean $q$-Poincar\'e quantum group in function algebra form}

For the construction of the $q$-Poincar\'e quantum group appropriate to
(\ref{eucmat}), we start with just the same formulae (\ref{vecpoi}) in
Section~3. For the construction there (from \cite{Ma:poi}) applies to any
braided covector space, and we have just seen in (\ref{eucbracovec}) that $\bar
A(R)$ can be put in this braided covector form. So the vectorial form of the
$q$-Poincar\'e quantum group for this case is just the same (\ref{vecpoi}) with
${\bf R}$ now determined from (\ref{eucbracovec}). The formula (\ref{veccoact})
is now
\eqn{eucveccoact}{\vecx\to \vecx\Lambda\dila+\vecp}
where $\vecp$ is the copy of algebra $\bar A(R)$ being used for momentum rather
than position $\vecx$.

Because the form of ${\bf R}$ for our particular example $\bar A(R)$ is simpler
than before, we have this time a simpler spinorial form. Firstly, we replace
$A({\bf R})$ by the tensor product Hopf algebra $A\tens A$ where each $A$ is a
quantum group obtained from $A(R)$, with generators $\vecs,\vect$ say. So the
spinorial spacetime rotation group is
\ceqn{eucspinlor}{ R\vecs_1\vecs_2=\vecs_2\vecs_1R,\quad \FRT,\quad
[\vect_1,\vecs_2]=0\\
\Delta\vecs=\vecs\tens\vecs,\quad \Delta\vect=\vect\tens\vect,\quad
\eps\vecs=\id=\eps\vect}
and further relations needed to have an antipode. The realisation of the
vectorial form in terms of the spinorial form is the same as (\ref{lambdast})
and the coaction takes the same form
\eqn{euccov}{\vecx\to \vecs^{-1}\vecx\vect\dila.}
For the standard example (\ref{Rsl2}) this fits with considerations for
$q$-Euclidean space in the appendix of \cite{CWSSW:ten}. As far as I know, the
general R-matrix setting using (\ref{eucmat}) is, however, due to
\cite{Ma:euc}. We take the tensor product $*$-structure
\eqn{eucspin*}{s^i{}_j{}^*=S^3s^j{}_i,\quad t^i{}_j{}^*=St^j{}_i,\quad
\dila^*=\dila}
again according to general theory in \cite{Ma:euc}. There one sees that an
extra automorphism $S^2$ in the definition of the $*$-structure in the first
copy of $A\tens A$ is needed for the coaction \cite[Eq. (11)]{Ma:euc} to be a
$*$-algebra map. In our case it becomes an extra $S^2$ on the $\vecs$ generator
for compatibility of (\ref{euccov}) with the Euclidean $*$-structure
(\ref{uni*}) on the space-time co-ordinates $\bar A(R)$, as well as for
$\Lambda^I{}_J{}^*=S\Lambda^J{}_I$. If we kept the original form (\ref{spin*})
with hermitian co-ordinates then such a system would be just our previous
$q$-Minkowski example in a twisted form.

The formulae (\ref{vecpoi}) then become the spinorial $q$-Poincar\'e quantum
group\cite{Ma:euc}
\ceqn{eucspinpoi}{ R_{21}\vecp_1\vecp_2=\vecp_2\vecp_1R,\quad
\vecp_1\vecs_2=\vecs_2\lambda^{-\h} R^{-1}\vecp_1,\quad
\vecp_1\vect_2=\lambda^\h \vect_2\vecp_1R\\
\vecp\dila=\lambda^{-1}\dila \vecp,\quad [\vecs,\dila]=[\vect,\dila]=0,\quad
\Delta\dila=\dila\tens\dila,\quad \eps\dila=1,\quad S\dila=\dila^{-1}\\
\Delta\vecp=\vecp\tens \vecs^{-1}(\ )\vect\dila+1\tens\vecp,\quad
\eps\vecp=0,\quad S\vecp=-\vecp\dila^{-1}\vect^{-1}(\ )\vecs}
where $\lambda$ is the quantum group normalisation constant of ${\bf R}$, which
is the square of that of $R$. Its value in the standard example is
$\lambda=q^{-1}$. As usual, one can derive both vectorial and spinorial forms
by the abstract bosonisation construction (\ref{cobos}) in Appendix~A. Finally,
the coaction on the spacetime generators in the spinorial form is of course
\eqn{eucspincoact}{ \vecx\to \vecs^{-1}\vecx\vect\dila+\vecp.}

\section{Euclidean $q$-Poincar\'e quantum group in enveloping algebra form}

For the enveloping algebra form for this $q$-Poincar\'e quantum group we use
(\ref{vecpoienv})  since this was a completely general construction for any
braided covector space\cite{Ma:poi}. We dualise each of the ingredients of the
semidirect product just as before. We use the pairing
$\<P^I,x_J\>=\delta^I{}_J$ between braided vectors and covectors in
(\ref{vecpair}), or $\<P_I,p_J\>=\eta_{IJ}$ in (\ref{lowmom}) for the lowered
index form of the Poincar\'e enveloping algebra (\ref{lowvecpoienv}). The only
difference from this part of Section~4 (up to and including (\ref{lowdif})) is
that the quantum metric tensor and ${\bf R}',{\bf R}$ now come from
(\ref{eucbracovec}) in Section~5,  where we cast $\bar A(R)$ as a braided
covector space. The action of the vectorial $q$-Poincar\'e enveloping algebra
on the spacetime co-ordinates takes the same form
\eqn{eucvecact}{\vecL^+_1\la\vecx_2=\vecx_2\lambda{\bf R}_{21},\quad
\vecL^-_1\la\vecx_2=\vecx_2\lambda^{-1}{\bf
R}^{-1},\quad\lambda^\xi\la\vecx=\lambda\vecx,\quad P^I\la x_J=\delta^I{}_J}
as before and necessarily makes the spacetime co-ordinates $\vecx$ into a
module algebra under it. The action of $P^I$ is by the braided differentials
$\bar\del^I$ as before defined with $[m;{\bf R}^{-1}_{21}]$ and obeying the
braided Leibniz rule (\ref{leib}) with respect now to $x_I$. Equivalently, when
there is a quantum metric $\eta$ the action of the lowered $P_I$ is by lowered
differentials $\bar\del_I$ obeying the braided covector algebra and Leibniz
rule
\ceqn{euclowleib}{ \bar\del_1\vecx_2-\lambda^2\vecx_2\bar\del_1{\bf R}=\eta}
as before. Indeed, everything for the vectorial $q$-Poincar\'e enveloping
algebras in Section~4 was for a general braided covector space.

We then use the specific form of ${\bf R}$ to give the spinorial description in
terms of $\vecx$ as a quantum matrix $\bar A(R)$. As in the preceding section,
it looks simpler than the braided matrix case in Section~4 because of the
simpler form of ${\bf R}$. This time the spinorial form of the spacetime
rotation enveloping
algebra is the tensor product Hopf algebra $H\tens H$ where each $H$ is dual to
$A$. So we take the two copies with FRT generators $\vecl^\pm$ and $\vecm^\pm$
as in (\ref{spinlorenv}) but now the matrix coproducts
\eqn{eucspinlorenv}{\Delta\vecl^\pm=\vecl^\pm\tens\vecl^\pm,\quad
\Delta\vecm^\pm=\vecm^\pm\tens\vecm^\pm,\quad \eps\vecl^\pm=\id,\quad
\eps\vecm^\pm=\id.}
The $*$-structure for the Euclidean picture dual to (\ref{eucspin*}) is
\eqn{eucspinenv*}{ l^\pm{}^i{}_j{}^*=S^{-1}l^\mp{}^j{}_i,\quad
m^\pm{}^i{}_j{}^*=Sm^\mp{}^j{}_i.}
The realisation of the vectorial spacetime rotation generators in the spinorial
ones is
\eqn{eucLlm}{
L^\pm{}^I{}_J=(S^{-1}l^\pm{}^{j_0}{}_{i_0})m^\pm{}^{i_1}{}_{j_1}.}
The spacetime-coordinates become a module algebra under the spacetime rotation
generators just as in (\ref{minkact}) by\cite{Ma:euc}
\eqn{eucact}{\vecl^+_1\la \vecx_2=\lambda^{-\h}R^{-1}_{21}\vecx_2,\quad
\vecl^-_1\la\vecx_2=\lambda^\h R\vecx_2,\quad \vecm^+_1\la \vecx_2=\vecx_2
\lambda^\h R_{21},\quad \vecm^-_1\la\vecx_2=\vecx_2 \lambda^{-\h}R^{-1}}
which is also the action we use on the lowered momentum generators
$P_I=P^{i_0}{}_{i_1}$ regarded now as a  matrix $\bar A(R)$. We add a dilaton
$\xi$ as before. Semidirect product by this action (and the induced semidirect
coproduct) according to the bosonisation construction in Appendix~A, or just
working from (\ref{lowvecpoienv}) for our particular ${\bf R}$ in
(\ref{eucbracovec}), immediately gives the spinorial form of the $q$-Poincar\'e
enveloping algebra as
\ceqn{eucspinpoienv}{R_{21}\vecP_1\vecP_2=\vecP_2\vecP_1R,\quad \vecl^+_1
\vecP_2=\lambda^{-\h}R^{-1}_{21}\vecP_2\vecl^+_1,\quad
\vecl^-_1\vecP_2=\lambda^\h R\vecP_2\vecl^-_1\\
 \vecm^+_1 \vecP_2=\vecP_2\lambda^\h R_{21}\vecm^+_1,\quad
\vecm^-_1\vecP_2=\vecP_2 \lambda^{-\h}R^{-1}\vecm^-_1,\quad
\lambda^\xi\vecP=\lambda^{-1}\vecP\lambda^\xi\\
{}[\xi,\vecl^\pm]=[\xi,\vecm^\pm]=0,\quad \Delta\xi=\xi\tens 1+1\tens\xi,\quad
\eps\xi=0,\quad S\xi=-\xi\\
\Delta \vecP=\vecP\tens 1+\lambda^\xi (\vecl^-(\ )S\vecm^-)\tens\vecP,\quad
\eps \vecP=0,\quad S\vecP=-\lambda^{-\xi}S(\vecl^-(\ )S\vecm^-)\vecP}
where $\vecl^-(\ )S\vecm^-$ has a space for the matrix indices of $\vecP$ to be
inserted. As far as I know, these are new formulae presented here for the first
time, although following directly from \cite{Ma:poi} as above, or by
dualisation of \cite{Ma:euc}. Finally, our $q$-Euclidean space $\bar A(R)$
becomes a module algebra under this $q$-Poincar\'e enveloping algebra by
\ceqn{eucspinact}{ \vecl^+_1\la \vecx_2=\lambda^{-\h}R^{-1}_{21}\vecx_2,\quad
\vecl^-_1\la\vecx_2=\lambda^\h R\vecx_2,\quad \vecm^+_1\la \vecx_2=\vecx_2
\lambda^\h R_{21},\quad \vecm^-_1\la\vecx_2=\vecx_2 \lambda^{-\h}R^{-1}\\
\lambda^\xi\la\vecx=\lambda\vecx,\quad  \vecP_1\la\vecx_2=\eta.}
The action of the lowered $P_I$ is by the lowered $\bar\del_I$ automatically
obeying relations and braided-Leibniz rule
\eqn{spinbardif}{R_{21}\bar\del_1\bar\del_2=\bar\del_2\bar\del_1R,\quad
\bar\del_1\vecx_2-\lambda^2R\vecx_2\bar\del_1R=\eta}
This is not a new equation but just (\ref{euclowleib}) in the matrix notation
according to the correspondence between the notations in (\ref{baraddstat}) and
(\ref{eucbracovec}).

The constructions in this spinorial setting work for general $q$-Hecke $R$,
with $\lambda$ the square of its quantum group normalisation constant. Finally,
let us note that there is an isomorphism of algebras between this
$q$-Poincar\'e enveloping algebra and the braided matrix case in Section~4 by
mapping $\vecP$ to $\vecP\vecl^-$. The two quantum groups are related by the
`quantum wick rotation' in \cite{Ma:euc}, namely by twisting
(cf. \cite{Dri:qua}) via the same quantum cocycle $\CR^{-1}$ which relates the
spacetime rotation quantum groups, but viewed now in the Euclidean Poincar\'e
enveloping algebra.

\section{Concluding remarks: $*$-structure and the dilaton problem}

We conclude here with two unconnected observations concerning problems of
current interest. Firstly, we know that our spacetime co-ordinates and Lorentz
algebras in all the above sections have reasonable $*$-structures and that the
Lorentz transformations preserve them. Hence it is natural to ask if the
Poincar\'e quantum group also has a natural $*$-structure. If we use the usual
axioms of a Hopf $*$-algebra then the answer appears in general to be no.

It does seem that one needs new axioms, albeit reducing to the usual ones when
$q=1$. A first step to formulating the correct axioms is in \cite{Ma:star}
where we study systematically $*$-structures on braided covector spaces. We
confirmed the axioms (\ref{delta*}) introduced in \cite{Ma:mec} and classify
the situations when they arise in the linear setting. We also explained there
that the construction of Poincar\'e groups by bosonisation would then inherit
natural $*$ properties with axioms to be elaborated elsewhere. Between the
release of  \cite{Ma:star} and the present paper, there appeared an
interesting preprint by Fiore\cite{Fio:man} in which the one-dimensional case
of a Poincar\'e algebra with two coproducts connected by $*$ was
considered directly.

In fact, the question of when bosonisations have natural structures as
$*$-algebras was already studied in \cite{Ma:mec} in order to have an
interpretation of the quantum double as quantum mechanics. For the enveloping
algebra Poincar\'e quantum group to be a $*$-algebra we need that the action
$\la$ of the `Lorentz' quantum group $H$ on the `momentum' braided group $C$ is
unitary in the standard Hopf algebra sense $(h\la c)^*=(Sh)^*\la c^*$ for $h$
in the Lorentz sector (including the dilaton) and $c$ in the momentum sector.
This {\em is} the case in all our examples above, for it just corresponds to
the coaction of the function algebra Lorentz quantum group being a $*$-algebra
homomorphism. The coalgebra however, is a semidirect product by the action
induced by the universal R-matrix of $H$ (see Appendix~A). We assume the latter
 is of real-type in the sense $\CR^{*\tens *}=\CR_{21}$, which again is the
case for our examples when $R$ is of real-type in the sense
$R^{\dagger\tens\dagger}=R_{21}$. This is true for (\ref{Rsl2}) when $q$ is
real.
One can easily see that in general the bosonisation $C\lbiprod H$ in this case
will {\em not} obey the usual axioms of a Hopf $*$-algebra. Instead we find

\begin{propos} In the general bosonisation theory\cite{Ma:bos} as in
Theorem~A.2 in the appendix, if the action of $H$ is `unitary' and its
universal R-matrix real-type as explained above then the coproduct and antipode
of $C\lbiprod H$ obey
\eqn{poi*}{ (*\tens*)\circ\Delta\circ *=\CR^{-1}(\tau\circ\Delta\ \, )\CR,\quad
\overline{\eps(\ )}=\eps\circ*,\quad  *\circ S\circ *=\cu^{-1}(S\ )\cu}
where $\cu=(S\CR\ut)\CR\uo$ and $\CR$ are viewed in the bigger algebra. We
propose to call a $*$-algebra with such an $\CR$ a {\em quasi-$*$} Hopf
algebra.
\end{propos}
The proof is easy by Hopf algebra techniques. Thus
\[(*\tens *)\circ\Delta c=\CR\ut^*c\Bo{}^*\tens (\CR\uo\la c\Bt)^*=\CR\umo
c\Bo{}^*\tens\CR\umt\la(c\Bt{}^*)=\bar\Delta(c^*)\]
where $\bar\Delta$ is the second `conjugate bosonisation' coproduct
(\ref{conjbos}), the second equality is our reality and unitarity assumption
and the third
is the $*$-axiom (\ref{delta*}) for braided groups. Here
$c\Bo\tens c\Bt$ is the braided coproduct of $C$ and $\CR\uo\tens\CR\ut$ the
universal R-matrix of $H$ (summations implicit). On the other hand,
$\bar\Delta$ is always twisiting equivalent to $\tau\circ\Delta$ by cocycle
$\CR^{-1}$. Since $H$ is a sub-Hopf algebra, we also have
$(\Delta\tens\id)\CR=\CR_{13}\CR_{23}$ and
$(\id\tens\Delta)\CR=\CR_{13}\CR_{12}$ when viewed in our bigger Hopf algebra.
We include these too in our characterisation of $C\lbiprod H$. So a quasi-$*$
Hopf algebra is like a quasitriangular Hopf algebra except that in one
axiom we replace $\Delta$ by $(*\tens *)\circ \Delta\circ *$.

The second coproduct $\bar \Delta$
is given for free in our braided approach as an automatic feature of the
theory: as explained in \cite{Ma:introp} every braided construction has a
conjugate one where we reverse the braid crossings. In the present case
$\Delta$ and $\bar\Delta$ coincide on the quantum group part $H$, where they
are both its usual coproduct. But on the braided group part $C$ they are more
like opposite (transposed) coproducts and indeed become that when $\CR=1$. This
is how we interpolate between the axioms (\ref{delta*}) for a braided group
(with a transposition $\tau$) and for a usual quantum group (without $\tau$)!

How does this second coproduct look for our $q$-Poincar\'e enveloping algebras
above? For the vectorial setting  (\ref{vecpoienv}) with upper indices $P^I$ it
is
\eqn{conjvec}{ \bar\Delta \vecP=\vecP\tens 1+\lambda^{-\xi}\vecL^+\tens
\vecP,\quad S\vecP=-\lambda^\xi(S\vecL^+)\vecP.}
If we use lowered covector indices $P_I$ then (\ref{lowvecpoienv}) becomes
\eqn{conjlowvec}{\bar\Delta P_I=P_I\tens1+\lambda^{-\xi}SL^+{}^A{}_I\tens
P_A,\quad SP_I=-\lambda^\xi(S^2 L^+{}^A{}_I)P_A.}
One can confirm (\ref{poi*}) directly for ${\bf R}$ of the appropriate reality
type and $*$ on the braided covectors as in our paper \cite{Ma:star}. For
example, in the Euclidean case it is $P_I^*=P^I$ while the Lorentz generators
obey $L^\pm{}^I{}_J{}^*=SL^\mp{}^J{}_I$. For the Minkowski $*$ in vectorial
form see \cite{Mey:new}. This all works generally for any
braided covector space with the correct reality properties.

For the Minkowski spinorial Poincar\'e enveloping algebra (\ref{spinpoienv}) in
Section~4, where $\vecP$ is a braided matrix $P^i{}_j$ in $B(R)$, the conjugate
quantum group structure is
\eqn{conjspin}{\bar\Delta \vecP=\vecP\tens 1+\lambda^{-\xi}\vecl^+\vecm^+(\
)(S_0\vecm^+)(S_0\vecl^-)\tens\vecP,\quad \bar S\vecP=-\lambda^\xi
S_0(\vecm^+\vecl^+(\ )(S_0\vecl^+)(S_0\vecm^+))\vecP}
where the space is for the matrix indices of $\vecP$ to be inserted and $S_0$
is the usual `matrix inverse' antipode. For the Euclidean one in
(\ref{eucspinpoienv}) in Section~7, where $P^i{}_j$ is in $\bar A(R)$, the
conjugate structure is
\eqn{conjeucspin}{\bar\Delta \vecP=\vecP\tens 1+\lambda^{-\xi}\vecl^+(\
)S\vecm^+\tens\vecP,\quad \bar S\vecP=-\lambda^\xi S(\vecl^+(\
)S\vecm^+)\vecP.}

Thus, the astute reader who was wondering why only $\vecL^-$  appeared in the
coproduct of $\vecP$ in Section~4, etc., in the construction of \cite{Ma:poi},
sees now that the symmetry is restored with the corresponding $\vecL^+$
appearing in the conjugate Poincar\'e quantum  group. This is why they are
connected by $*$ as in (\ref{poi*}). This suggests that these are indeed very
reasonable axioms for our setting. Finally, while the Poincar\'e generators
before were represented on spacetime co-ordinates by means of $\bar\del$
braided differentials, it is easy to see that with the conjugate coproduct the
same algebra acts covariantly with the usual braided differentials $\del$ from
\cite{Ma:fre}. We saw already in \cite[Sec. V]{Ma:fre} that if $V({\bf R}',{\bf
R})$ is a braided vector space of differentials (our momentum sector) then its
`conjugate' or opposite braided group is $V({\bf R}',{\bf R}_{21}^{-1})$
obeying the same algebra with reversed braiding. When we bosonise, it means
that our one $q$-Poincar\'e enveloping algebra extends to products of spacetime
generators in {\em two} different ways (\ref{leib}), one with the Leibniz rule
for $\bar\del$ and the other with the conjugate-braiding Leibniz rule for
$\del$. The two coincide in the triangular or unbraided case so the distinction
is not visible classically. In physical terms it means that $q$-deformed
geometry in the braided approach is naturally `split' into two geometries
related by braid-crossing reversal symmetry in the constructions.

It is clear that Proposition~8.1 also solves in principle the important problem
of how to tensor product `unitary' representations of the $q$-Poincar\'e group,
which we need for physics. That the conjugated coproduct $(*\tens
*)\circ\Delta\circ *$ is twisting equivalent to $\tau\circ\Delta$ means that
the tensor product of two representations which are unitary in the sense
$\rho(h^*)^\dagger=\rho(\sigma h)$ (where $\sigma$ is a fixed algebra and
anticoalgebra map) will remain so, {\em but up to isomorphism}. In our case
the isomorphism is given by the action of $\CR^{-1}$ and is a new
physical effect which is absent when $\CR=1$. For example, we can define a {\em
braided-unitary} representation to be a pair consisting of $V$ on which our
Poincar\'e quantum group acts and a semilinear form $(\ ,\ )_V$ on it such that
$(h^*\la v,v')_V=(v,(\sigma h)\la v')_V$ for all $v,v'\in V$ and $h$ in our
Poincar\'e (or other quasi-$*$ Hopf) algebra. Then one can see that
\eqn{tensuni}{V\tens W,\quad (v\tens w,v'\tens w')_{V\tens W}=(\CR\umt\la
v,v')_V(\CR\umo\la w,w')_W}
is again a braided-unitary representation, where we act on tensor products in
the usual way via the coproduct $\Delta$. The definition is associative using
the other properties of $\CR$. We can also eliminate $\sigma$ by taking a more
categorical line with one input of $(\ ,\ )$ living in the category with
opposite tensor product. We do not, however, assume that the semilinear
forms are conjugate symmetric (like a usual Hilbert space inner product) since
this is not in general preserved by (\ref{tensuni}). The problem of
constructing and perhaps classifying such braided-unitary representations for
our particular Poincar\'e algebras will be addressed elsewhere.

For our second topic we recall that the appearance of the dilaton $\dila$ or
$\xi$ in the $q$-Poincar\'e quantum group is an unexpected feature noted
already in \cite{SWW:inh}\cite{OSWZ:def}, and explained in terms of the quantum
group normalisation constant in the general construction \cite{Ma:poi}. We
observe  now  that, while not solving anything, we can systematically remove it
from all our q-Poincar\'e quantum group formulae
in our constructions above if we pay a certain price. The price is that we
obtain then not a quantum group in the usual sense but a $\Z$-graded or
$\C$-statistical braided group\cite{Ma:csta}. This is a braided group of the
kind where the statistics are just a power of a generic factor, say $\lambda$.
This class includes as special cases superquantum groups and anyonic quantum
groups\cite{Ma:any}\cite{MaPla:any}. If $B$ is such a braided group, then its
bosonisation consists of adjoining  a new group-like generator $g$, say, with
commutation relations which remember the grading $|b|\in\Z$. It has
structure\cite{Ma:bos}\cite{Ma:introp}
\ceqn{gbos}{gb=\lambda^{|b|}bg,\quad \Delta g=g\tens g,\quad\eps g=1,\quad
Sg=g^{-1}\\
\Delta b=b\Bo g^{|b\Bt|}\tens b\Bt,\quad Sb=g^{-|b|}\und S b}
where $b\Bo\tens b\Bt$ is the braided coproduct and $\und S$ the braided
antipode of $B$, and $|\ |$ is the degree of a
homogeneous element. We use once again the general theory in Appendix~A, either
as (a left-handed version of) the bosonisation (\ref{cobos}) or in terms of the
enveloping algebra version (\ref{bos}) with $H$ the quantum line.

Now, comparing this formula with the formulae for our $q$-Poincar\'e enveloping
algebras
(\ref{vecpoienv}),(\ref{lowvecpoienv}),(\ref{spinpoienv}),(\ref{eucspinpoienv})
we see that in each case they are of the above form with
\eqn{poigrad}{ g=\lambda^\xi,\quad |\vecL^\pm|=|\vecl^\pm|=|\vecm^\pm|=0,\quad
|\vecP|=-1;\quad |\vecu|=|\vecx|=1}
and with $B$ defined as the $\Z$-graded braided group given by the same
formulae with $\lambda^\xi$ omitted. So these $B$ are braided $q$-Poincar\'e
groups. They act on the spacetime in a way that preserves grading also, i.e.
these become $\Z$-graded module algebras with the grading as shown, etc.
Similarly in right-handed conventions for the function algebra quantum
Poincar\'e groups in Sections~3,6 with $\dila$ omitted. Everything works as in
supersymmetry except that the statistical factor $-1$ is replaced by
$\lambda$\cite{Ma:csta}. In physical terms, the role played mathematically in
supersymmetry by fermionic degree is played now by the physical scale
dimension.
This is another example of the unification of different physical concepts made
possible by quantum groups and braided groups.

\appendix
\section{Appendix: the abstract bosonisation theory}

Here we collect some basic formulae from the abstract theory of
bosonisation\cite{Ma:bos}\cite{Ma:mec}. Conceptually, bosonisation is a
generalisation of the Jordan-Wigner transform for turning fermionic systems
into bosonic ones. Recall that this is done by adjoining the degree or grading
operator. In braided geometry the role of the $\Z_2$ grading of supersymmetry
is played by the background quantum group symmetry, which in the above context
is the spacetime rotation quantum group (the concepts of supersymmetry and
Lorentz invariance are unified\cite{Ma:introp}).

We assume that the reader is familiar with the definition of a quasitriangular
Hopf algebra $H$ in \cite{Dri} with `universal R-matrix'
$\CR=\CR\uo\tens\CR\ut$ in $H\tens H$ (summation of terms implicit), and the
dual notion of a dual-quasitriangular Hopf algebra $A$ with
dual-quasitriangular structure or `universal R-matrix functional' $\CR:A\tens
A\to \C$ in \cite{Ma:eul}\cite{Ma:bg}. The latter is characterised by the
axioms
\ceqn{dqua}{ \CR(a,bc)=\CR(a\o,c)\CR(a\t,b),\quad
\CR(ab,c)=\CR(a,c\o)\CR(b,c\t)\\ b\o a\o \CR(a\t,b\t)=\CR(a\o,b\o)a\t b\t.}
The coproducts are denoted $\Delta a=a\o\tens a\t$, etc. (summation implicit).
We also assume that the reader is familiar with the basic notion of a braided
group $B$ or braided-Hopf algebra\cite{Ma:bra}\cite{Ma:bg}. Introductions are
in \cite{Ma:introp}\cite{Ma:introm}\cite{Ma:varen}.

The first theorem is that if $B$ is a braided group living in the braided
category of (say right) $A$-comodules (i.e. an object which is totally
covariant under a coaction of $A$) then there is an ordinary Hopf algebra
$A\rbiprod B$ constructed as follows\cite{Ma:mec}. As a coalgebra we make a
semidirect coproduct by the right coaction of $A$. We also use the universal
R-matrix functional to turn the right coaction into a right action of $A$ by
evaluation in one input, and make an algebra semidirect product by this action.
In concrete terms:

\begin{theorem}\cite{Ma:mec} cf.\cite{Ma:bos} The algebra $A\rbiprod B$
generated by $B,A$ with cross relations, coproduct and antipode
\eqn{cobos}{ ba=a\o b\bo \CR(b\bt,a\t),\quad \Delta b=b\Bo\bo\tens b\Bo\bt
b\Bt,\quad S b=(\und S b\bo)Sb\bt}
for all $a\in A, b\in B$, is a Hopf algebra.
\end{theorem}
Here $\und\Delta b=b\Bo\tens b\Bt$ is the braided coproduct of $B$, $\und S$
its braided antipode and $b\bo\tens b\bt\in B\tens A$ denotes the output of the
coaction of $A$. The coproduct and antipode of $A$ are not modified (so $A$ is
a sub-Hopf algebra).

The second theorem (actually proven first in \cite{Ma:bos} with the above
easily obtained as dual to it) is that if $C$ is some braided group living in
the braided category of (say, left) $H$-modules (i.e. an object which is
totally covariant under an action of $H$) then there is an ordinary Hopf
algebra $C\lbiprod H$ constructed as follows\cite{Ma:bos}. As an algebra we
make a semidirect product by the left action of $H$. Also, we use the universal
R-matrix to turn the left action into a left coaction by letting its part
living in the first factor of $H$ act. In concrete terms:

\begin{theorem}\cite{Ma:bos} The algebra $C\lbiprod H$ generated by $H,C$ with
cross relations, coproduct and antipode
\ceqn{bos}{ hc=(h\o\la c)h\t,\quad \Delta c=c\Bo \CR\ut\tens \CR\uo\la
c\Bt,\quad Sc=(\cu \CR\uo\la\und Sc) S\CR\ut}
for all $h\in H, c\in C$, is a Hopf algebra.
\end{theorem}
Here $\cu=(S\CR\ut)\CR\uo$, $\und\Delta c=c\Bo\tens c\Bt$ the braided coproduct
of $C$, $\und S$ its braided antipode and $\la$ denotes the left action of $H$.
The coproduct and antipode of $H$ are not modified (so $H$ is a sub-Hopf
algebra).

It should be perfectly clear that these two constructions are conceptually dual
to one another. So if $C=B^\star$ (the braided group dual to $B$) and $H=A^*$
(the quantum group dual to $A$) then $C\lbiprod H=(A\rbiprod B)^*$ as usual
quantum groups.

\begin{propos}cf.\cite{Ma:mec} The duality pairing of $C\lbiprod H$ and
$A\rbiprod B$ between the various subalgebras is
\eqn{bos-cobos}{ \<c,a\>=\eps(c)\eps(a),\quad \<h,a\>={\rm usual},\quad
\<c,b\>=\ev(\und S^{-1}c,b),\quad \<h,b\>=\eps(h)\eps(b)}
where $\eps$ is the counit of the quantum group or braided group, `usual' means
the pairing between $H,A$ as usual quantum groups dual to each other, $\und S$
is the braided antipode (which we assume invertible) of $C$, and $\ev$ is the
braided-group duality evaluation pairing.
\end{propos}

Recall\cite{Ma:introp} that $\ev$ obeys slightly different axioms to a usual
quantum group pairing. Indeed, $\ev(\und S^{-1}(\ ),(\ ))$ obeys axioms more
like the usual axioms and reduces to them when the braiding is trivial. We have
slightly reworked \cite{Ma:mec} where the duality was given explicitly when
$B=C^\star$ rather than $C=B^\star$ as here. Both statements are true. Also
coming out of the bosonisation theory is a canonical coaction of $A\rbiprod B$
on $B$ or, by duality in the setting above, an action of $C\lbiprod H$ on $B$.
By its very definition in categorical terms, the ordinary (co)representations
of the bosonised Hopf algebra are in 1-1 correspondence with the braided
(co)-representations of the braided group before bosonisation\cite{Ma:bos}.
Obviously $B$ coacts on itself by its braided coproduct. So the general theory
gives at once:

\begin{cor}cf.\cite{Ma:poi} $B$ is a right $A\lbiprod B$-comodule algebra by
\eqn{boscoact}{b\to b\Bo\bo\tens b\Bo\bt b\Bt}
and in the setting above a left $C\lbiprod H$-module algebra by
\eqn{bosact}{h\la b=b\bo\<h,b\bt\>,\quad c\la b=b\Bo\bo\ev(\und
S^{-1}c,b\Bt\bo)\CR^{-1}(b\Bt\bt,b\Bo\bt)}
\end{cor}
Conceptually, the action of $H$ on $B$ is just the action corresponding to the
coaction of $A$ assumed when we said that $B$ was $A$-covariant to begin with.
The action of $C$ in abstract terms is
\eqn{Lreg}{ c\la b=(\ev(\und S^{-1}c,(\
))\tens\id)\circ\Psi^{-1}\circ\und\Delta b}
which has a braided picture when we write $\Psi$ as a braid crossing and
$\ev=\cup$. It is a left-handed version of the right coregular representation
${\rm Reg}^*$ introduced and studied in \cite{Ma:introp}\cite{Ma:lie}, and
always makes a braided group $B$ a braided module algebra under its dual
braided group.

The combination $\bar\Delta=\Psi^{-1}\circ\Delta$ is studied in
\cite{Ma:introp} as the {\em naive opposite coproduct}. It is naive because it
does not make the algebra of $C$ into a braided group in our original braided
category but rather into a braided group $\bar C$, say, living in the
`conjugate' braided category with inverse transposed
braiding\cite[Lemma~4.6]{Ma:introp}. $\und S^{-1}$ becomes its braided
antipode.  In concrete terms it means that the braided group $\bar C$ is no
longer properly covariant under $H$ (with the correct induced braiding) but
under the quantum group $H$ equipped with $\CR_{21}^{-1}$ instead for its
universal R-matrix. Let us denote the latter by $\bar H$. As a Hopf algebra it
coincides with $H$, but has `conjugate' $\CR$.
\begin{cor} Every bosonisation $C\lbiprod H$ has a second `conjugate' coproduct
and antipode on the same algebra,
\eqn{conjbos}{\bar\Delta c=\CR\umo c\Bt\tens \CR\umt\la c\Bo,\quad \bar
Sc=(\CR\ut\cv^{-1}\la \und S^{-1}c)\CR\uo}
where $\cv=\CR\uo S\CR\ut$.
\end{cor}
This is just the bosonisation $\bar C\lbiprod\bar H$ of $\bar C$ from
Theorem~A.2, written (using the algebra relations) in terms of $H,C$. It has
the same algebra as $C\lbiprod H$ but is generally a different Hopf algebra.
Likewise every bosonisation $A\rbiprod B$ has a conjugate $\bar A\rbiprod\bar
B$ with the same coalgebra and a different product. This kind of braid-crossing
reversal symmetry is an intrinsic feature of braided group
theory\cite{Ma:introp}. It is easy to further recognise these second coproducts
and products as twisiting equivalent (cf. \cite{Dri:qua}) to the
opposite coproduct or product respectively. Thus
$\bar\Delta=\CR^{-1}(\tau\circ\Delta\ )\CR$ by an elementary computation using
the algebra relations.

This summarises the relevant parts of the abstract theory of bosonisation of
braided groups into ordinary quantum groups. One of the first applications in
physics was to the construction of $q$-Poincar\'e quantum group function
algebras for any braided covector space $B=\Vhaj({\bf R}',{\bf
R})$\cite{Ma:poi}. The latter lives in the braided category of
$\widetilde{A}$-comodules by (\ref{veccov}), where $A$ is obtained from $A({\bf
R})$ with universal R-matrix functional in \cite{Ma:lin}\cite{Ma:qua} and we
extend by a dilaton with universal R-matrix functional
$\CR(\dila,\dila)=\lambda^{-1}$. We obtain from (\ref{cobos}) the formula
(\ref{vecpoi}) in Section~3 as the Hopf algebra $\widetilde{A}\rbiprod
\Vhaj({\bf R}',{\bf R})$; see \cite{Ma:poi}.

For the spinorial examples, we let instead $A$ be a quantum group obtained from
$A(R)$ and take $B(R)$ in the category of $\widetilde{A\dcross A}$-comodules by
(\ref{minkcov}) to obtain (\ref{spinpoi}) as $\widetilde{A\dcross A}\rbiprod
B(R)$; see \cite{MaMey:bra}. Likewise, we take $\bar A(R)$ in the category of
$\widetilde{A\tens A}$ comodules under (\ref{euccov}) to obtain
(\ref{eucspinpoi}) as $\widetilde{A\tens A}\rbiprod \bar A(R)$; see
\cite{Ma:euc}.

We also gave the dual construction in \cite{Ma:poi} to obtain Poincar\'e
quantum enveloping algebras dual by (\ref{bos-cobos}) to the above examples.
Thus $C=V({\bf R}',{\bf R})$ is the braided vector space, dual to the braided
covectors above using braided-differentiation\cite{Ma:fre}\cite{KemMa:alg}. The
symbol $\Vhaj$ denotes the predual, i.e. $V=(\Vhaj)^*$. The generators are dual
spaces to each other (as for usual vectors and covectors).
This braided group lives in the same category of $\widetilde{A}$ comodules, or
equivalently, $\widetilde {H}$-modules, where $\widetilde{H}$ is dual to
$\widetilde{A}$ and $A$ is from $A({\bf R})$. The dilaton contributes
$\lambda^{-\xi\tens\xi}$ to the universal R-matrix. We obtain from (\ref{bos})
the formula (\ref{vecpoienv}) in Section~4 as $V({\bf R}',{\bf
R})\lbiprod\widetilde{H}$; see \cite{Ma:poi}. We obtain from (\ref{bos-cobos})
the explicit duality pairing (\ref{vecpair}).

Likewise for the spinorial examples, we let $H$ be dual to $A$ obtained from
$A(R)$, and obtain (\ref{spinpoienv}) as the Hopf algebra
$B(R)^*\lbiprod\widetilde{H\codcross H}$ and (\ref{eucspinpoienv}) as the Hopf
algebra $\bar A(R)^*\lbiprod\widetilde{H\tens H}$. These are canonical
constructions independent of any quantum metric. When there is a quantum
metric, as in the explicit $q$-Minkowski and $q$-Euclidean examples, we have
$B(R)^*\isom B(R)$ and $\bar A(R)\isom \bar A(R)$ in which final form we wrote
these examples. They are equivalent to constructing the bosonisations
$B(R)\lbiprod\widetilde{H\codcross H}$ and $\bar A(R)\lbiprod \widetilde{H\tens
H}$ directly.

Finally, all these examples come equipped with a canonical (co)action from
Corollary~A.4 on the braided covectors regarded as spacetime co-ordinates. This
gives the explicit formulae (\ref{veccoact}) and (\ref{vecact}) in general and
(\ref{spinact}), (\ref{eucspinact}) on $B(R), \bar A(R)$. We systematically
dropped a minus sign coming from the braided antipode on $\vecP$ in
(\ref{Lreg}). In addition, all our $q$-Poincar\'e quantum groups automatically
come with conjugates from Corollary~A.5, as discussed in Section~8.

We have concentrated here on the applications of bosonisation to the
construction of $q$-Poincar\'e quantum groups. Other interesting applications
are in \cite{MacMa:str}\cite{MaPla:uni} and the theory of differential calculus
on quantum groups in \cite{SchZum:bra}\cite{Dra:bra}.


\end{document}